# Advanced Wigner Distribution and Ambiguity Function in the Quadratic-phase Fourier Transform Domain: Mathematical Foundations and Practical Applications

**Aamir H. Dar\* and Neeraj Kumar Sharma**

*Mehta Family School of Data Science & Artificial Intelligence*
*Indian Institute of Technology Guwahati, Guwahati-781039, India*
*E-mail:  aamir740@rnd.iitg.ac.in; neerajs@iitg.ac.in*

ABSTRACT. In non-stationary signal processing, prior work has incorporated the quadratic-phase Fourier transform (QPFT) into the ambiguity function (AF) and Wigner distribution (WD) to enhance their performance. This paper introduces an advanced quadratic-phase Wigner distribution and ambiguity function (AQWD/AQAF), extending classical WD/AF formulations. Key properties, including the Moyal formula, anti-derivative property, shift, conjugation symmetry, and marginal properties, are established. Furthermore, the proposed distributions demonstrate improved effectiveness in linear frequency-modulated (LFM) signal detection. Simulation results confirm that AQWD/AQAF outperforms both traditional WD/AF and existing QPFT-based WD/AF methods in detection accuracy and overall performance.

## 1. Introduction

The Fourier transform (FT) is a technique for analyzing the frequency spectrum of stationary signals. It has enabled analyze stationary signals in the frequency and time domains, independently. The FT has contributed to numerous developments in diverse fields of engineering and sciences. However, the FT is less effective in interpreting the temporally evolving spectrum in non-stationary signals. Castro et al. [1, 2] proposed the quadratic-phase Fourier transform (QPFT). This generalized version of the FT provides an approach for unified treatment of both stationary and non-stationary signals.

The QPFT of a signal $f(t)$, with a given set of real parameters $\Lambda = (A, B, C, D, E,$ and $B \neq 0)$ is defined as follows [3, 4]

$$\mathbb{Q}_\Lambda\{f(t)\}(\nu) = \frac{1}{\sqrt{2\pi}} \int_\mathbb{R} f(t) \mathcal{K}_\Lambda(\nu, t) \mathrm{d}t, \qquad (1.1)$$

where the QPFT kernel $\mathcal{K}_\Lambda(\nu, t)$ is given by

$$\mathcal{K}_\Lambda(\nu, t) := \sqrt{\frac{B}{i}} \mathrm{e}^{i\left(A\nu^2 + Bt\nu + Ct^2 + D\nu + Et\right)}. \qquad (1.2)$$

Many well-known integral transforms, such as the FT, linear canonical transform (LCT), fractional Fourier transform (FrFT), offset linear canonical transform (OLCT), and so on, are embodied in the QPFT. The QPFT has demonstrated its reliability as a tool for the efficient representation of stationary and non-stationary signals requiring multiple controllable parameters in various engineering and scientific fields, such as harmonic analysis, image processing, sampling, kernel replicating theory, and several others [5]-[8].

In contrast to frequency domain representations, there are also approaches for obtaining time-frequency distributions which benefit the analysis of non-stationary signals [9]. The classical Wigner distribution (WD) and the classical ambiguity function (AF) are among the fundamental nonparametric time-frequency analytic tools. They are mainly used in applications to identify LFM signals to examine the time-frequency characteristics of non-stationary signals [10]-[17]. The WD and the AF of $f \in L^2(\mathbb{R})$ are defined by [18, 19, 20]

$$\mathcal{W}_f(t,\nu) = \int_{\mathbb{R}} f\left(t+\frac{\Upsilon}{2}\right) f^*\left(t-\frac{\Upsilon}{2}\right) e^{-i\nu\Upsilon} d\Upsilon, \tag{1.3}$$

$$\mathcal{A}_f(\Upsilon,\nu) = \int_{\mathbb{R}} f\left(t+\frac{\Upsilon}{2}\right) f^*\left(t-\frac{\Upsilon}{2}\right) e^{-i\nu t} dt, \tag{1.4}$$

where the superscript "*" denotes the complex conjugation. Moreover, with the help of classical convolution defined as:

$$(f * g)(t) = \int_{\mathbb{R}} f(\Upsilon) g(t-\Upsilon) d\Upsilon, \quad f, g \in L^2(\mathbb{R}), \tag{1.5}$$

equations (1.3) and (1.4) can be expressed as follows:

$$\mathcal{W}_f\left(\frac{t}{2},\nu\right) = 2[f(t)e^{-i\nu t}] * [f^*(t)e^{i\nu t}], \tag{1.6}$$

$$\mathcal{A}_f(\Upsilon, 2\nu) = [f(\Upsilon)e^{-i\nu\Upsilon}] * [f^*(-\Upsilon)e^{i\nu\Upsilon}]. \tag{1.7}$$

Although the quadratic phase Fourier transform (QPFT) is effective in identifying linear frequency-modulated (LFM) signals due to its energy accumulation properties in specific QPFT domains, the combined use of QPFT with Wigner distribution (WD) and ambiguity function (AF) further enhances LFM signal detection. The initial exploration of QPFT-based WD/AF, adhering to the classical framework of WD and AF, was carried out by the authors in [21, 22]. More recently, Bhat and Dar [23] introduced a variant, the quadratic-phase Wigner distribution and ambiguity function (QWD/QAF). The WD and AF associated with QPFT, of a finite energy signal $f \in L^2(\mathbb{R})$, is expressed as,

$$\text{QWD}_f^{\Lambda}(t,\nu) = \sqrt{\frac{B}{2\pi i}} \int_{\mathbb{R}} f\left(t+\frac{\Upsilon}{2}\right) f^*\left(t-\frac{\Upsilon}{2}\right) e^{i(A\Upsilon^2 + B\nu\Upsilon + C\nu^2 + D\Upsilon + E\nu)} d\Upsilon, \tag{1.8}$$

$$\text{QAF}_f^{\Lambda}(\Upsilon,\nu) = \sqrt{\frac{B}{2\pi i}} \int_{\mathbb{R}} f\left(t+\frac{\Upsilon}{2}\right) f^*\left(t-\frac{\Upsilon}{2}\right) e^{i(At^2 + B\nu t + C\nu^2 + Dt + E\nu)} dt. \tag{1.9}$$

Further generalization of WD/AF using the QPFT domain and its applications is explored in [24, 25, 26]. Using the flexibility of the quadratic-phase kernel and conventional convolution operator, this paper provides an alternative description of WD and AF associated with QPFT. Subsequently, it presents interesting properties and potential applications compared to [21]-[26]. The following are the key contributions of this paper.

- Introduces a novel Wigner distribution (WD) and ambiguity function (AF) in the quadratic-phase Fourier transform (QPFT) domain, referred to as the advanced quadratic-phase Wigner distribution and ambiguity function (AQWD/AQAF).
- Conducts analysis of the fundamental properties of AQWD/AQAF, including shift, conjugation symmetry, marginal properties, the Moyal formula, and the anti-derivative property.
- Demonstrates the applicability of AQWD/AQAF in identifying linear frequency-modulated (LFM) signals, both single and multi-component, and provides simulation results confirming its superior performance over traditional WD/AF and existing QPFT-based QWD/QAF methods.

The rest of the paper is organized as follows. Section 2 introduces the definition of AQWD/AQAF, establishes its connection to classical WD and AF, and explores special cases of the proposed formulation. Section 3 presents a detailed discussion of key properties, including shift, conjugation symmetry, marginal properties, the Moyal formula, and the anti-derivative property. The implications of AQWD/AQAF for detecting single and multi-component LFM signals are examined in Section 4, supported by simulation results that validate the proposed methodology. Finally, Section 5 concludes the paper by summarizing key findings and outlining potential directions for future research.

## 2. Advanced Quadratic-phase Wigner Distribution and Ambiguity Function

This section introduces the proposed AQWD and AQAF, and their special cases.

### 2.1. Definition of the AQWD and AQAF.

Thanks to conventional convolution (1.5), expressions of classical WD and AF can be modified as:

When $e^{-i\nu t}$ is substituted with $\mathcal{K}_\Lambda(\nu, t)$ and $e^{i\nu t}$ with $\mathcal{K}_\Lambda^*(t, \nu)$ in (1.6), and the variable $z = \frac{t}{2} + \frac{\Upsilon}{2}$ is changed, we get

$$2\left[f(t)\mathcal{K}_\Lambda(\nu, t)\right] * \left[f(t)\mathcal{K}_A(t, \nu)\right]^*$$
$$= 2\int_{\mathbb{R}} f(z) f^*(t-z) \mathcal{K}_\Lambda(\nu, z) \mathcal{K}_\Lambda^*(t-z, \nu) \, dz$$
$$= |B| \int_{\mathbb{R}} f\left(\frac{t}{2} + \frac{\Upsilon}{2}\right) f^*\left(\frac{t}{2} - \frac{\Upsilon}{2}\right) \mathcal{K}_\Lambda\left(\nu, \frac{t}{2} + \frac{\Upsilon}{2}\right) \mathcal{K}_\Lambda^*\left(\frac{t}{2} - \frac{\Upsilon}{2}, \nu\right) d\Upsilon$$
$$= |B|e^{i\left[(A-C)\nu^2 + (D-E)\nu\right]} \int_{\mathbb{R}} f\left(\frac{t}{2} + \frac{\Upsilon}{2}\right) e^{i\left[C\left(\frac{t}{2} + \frac{\Upsilon}{2}\right)^2 + \left(\frac{t}{2} + \frac{\Upsilon}{2}\right)E\right]}$$
$$\times f^*\left(\frac{t}{2} - \frac{\Upsilon}{2}\right) e^{-i\left[A\left(\frac{t}{2} - \frac{\Upsilon}{2}\right)^2 + D\left(\frac{t}{2} - \frac{\Upsilon}{2}\right)\right]} e^{iB\nu\Upsilon} d\Upsilon.$$
(2.1)

Similarly, by replacing $e^{-i\nu\Upsilon}$ with $\mathcal{K}_\Lambda(\nu, \Upsilon)$ and $e^{i\nu\Upsilon}$ with $\mathcal{K}_\Lambda^*(\Upsilon, \nu)$ in (1.7) and the variable $z_1 = t + \frac{\Upsilon}{2}$ is changed, we get

$$[f(\Upsilon)\mathcal{K}_\Lambda(\nu, \Upsilon)] * [f^*(-\Upsilon)\mathcal{K}_\Lambda^*(\Upsilon, \nu)]$$
$$= \int_{\mathbb{R}} f(y)\mathcal{K}_\Lambda(\nu, y) f^*(z_1 - \Upsilon)\mathcal{K}_\Lambda^*(\Upsilon - z_1, \nu) \, dz_1$$
$$= |B| \int_{\mathbb{R}} f\left(t + \frac{\Upsilon}{2}\right) f^*\left(t - \frac{\Upsilon}{2}\right) \mathcal{K}_\Lambda\left(\nu, t + \frac{\Upsilon}{2}\right) \mathcal{K}_\Lambda^*\left(\frac{\Upsilon}{2} - t, \nu\right) dt$$
$$= |B|e^{i\left[(A-C)\nu^2 + (D-E)\nu\right]} \int_{\mathbb{R}} f\left(t + \frac{\Upsilon}{2}\right) e^{i\left[C\left(t + \frac{\Upsilon}{2}\right)^2 + \left(t + \frac{\Upsilon}{2}\right)E\right]}$$
$$\times f^*\left(t - \frac{\Upsilon}{2}\right) e^{-i\left[A\left(t - \frac{\Upsilon}{2}\right)^2 + D\left(t - \frac{\Upsilon}{2}\right)\right]} e^{i2B\nu t} dt.$$
(2.2)

The following definition of advanced WD and AF related to the QPFT (AQWD/AQAF) is derived from the aforementioned equations (2.1) and (2.2).

**Definition 2.1.** For a real parametric Set $\Lambda = (A, B, C, D, E), B \neq 0$, the AQWD and AQAF of a signal $f \in L^2(\mathbb{R})$ are defined as

$$\mathcal{QW}_f^\Lambda(t,\nu) = |B|e^{i[(A-C)\nu^2+(D-E)\nu]}\int_\mathbb{R} f_{C,E}\left(t+\frac{\Upsilon}{2}\right)f^*_{A,D}\left(t-\frac{\Upsilon}{2}\right)e^{iB\nu\Upsilon}d\Upsilon, \quad (2.3)$$

$$\mathcal{QA}_f^\Lambda(\Upsilon,\nu) = |B|e^{i[(A-C)\nu^2+(D-E)\nu]}\int_\mathbb{R} f_{C,E}\left(t+\frac{\Upsilon}{2}\right)f^*_{A,D}\left(t-\frac{\Upsilon}{2}\right)e^{iB\nu t}dt \quad (2.4)$$

where $f_{p,q}(t) := f(t)e^{i(pt^2+qt)}$.

Furthermore, it may be stated that the AQWD and AQAF are scaling and modulated variants of the cross-term WD and AF in the time and frequency variables, respectively via following formula

$$\begin{aligned}
\mathcal{QW}_f^\Lambda(t,\nu) &= |B|e^{i[(A-C)\nu^2+(D-E)\nu]}\mathcal{W}_{f_{C,E},f_{A,D}}(t,-B\nu)\\
&= |B|e^{i[(A-C)\nu^2+(D-E)\nu]}\mathcal{W}_{f,f_{A-C,D-E}}(t,-(E+2Ct+B\nu)),
\end{aligned}$$

$$\begin{aligned}
\mathcal{QA}_f^\Lambda(\Upsilon,\nu) &= |B|e^{i[(A-C)\nu^2+(D-E)\nu]}\mathcal{W}_{f_{C,E},f_{A,D}}(\Upsilon,-B\nu)\\
&= |B|e^{i[(A-C)\nu^2+(D-E)\nu]}\mathcal{W}_{f,f_{A-C,D-E}}(t,-(E+2C\Upsilon+B\nu)).
\end{aligned}$$

The following well known integral transforms are obtained when parameter set $\Lambda = (A, B, C, D, E)$ has certain particular forms:

- By choosing the parametric set $\Lambda = (\frac{D}{2B}, \frac{-1}{B}, \frac{A}{2B}, 0, 0)$, the AQWD and AQAF simply yield novel Wigner distribution and ambiguity function in the LCT domain by L.T. Minh [27]

$$\mathcal{QW}_f^\Lambda(t,\nu) = \frac{1}{|B|}e^{\frac{i(D-A)}{2B}\nu^2}\int_\mathbb{R} f_{\frac{A}{2B}}\left(t+\frac{\Upsilon}{2}\right)f^*_{\frac{D}{2B}}\left(t-\frac{\Upsilon}{2}\right)e^{-\frac{i}{B}\nu\Upsilon}d\Upsilon, \quad (2.5)$$

$$\mathcal{QA}_f^\Lambda(\Upsilon,\nu) = \frac{1}{|B|}e^{\frac{i(D-A)}{2B}\nu^2}\int_\mathbb{R} f_{\frac{A}{2B}}\left(t+\frac{\Upsilon}{2}\right)f^*_{\frac{D}{2B}}\left(t-\frac{\Upsilon}{2}\right)e^{-\frac{i}{b}\nu t}dt \quad (2.6)$$

where $f_p(t) := f(t)e^{ipt^2}$.

- When $A = C$ and $D = E$ in $\Lambda = (A, B, C, D, E)$, the proposed AQWD and AQAF yields recently introduced WD and AF in QPFT domain ( see [26] ) as

$$\mathcal{QW}_f^\Lambda(t,\nu) = |B|\int_\mathbb{R} f\left(t+\frac{\Upsilon}{2}\right)f^*\left(t-\frac{\Upsilon}{2}\right)e^{i(2At+B\nu+D)\Upsilon}d\Upsilon, \quad (2.7)$$

$$\mathcal{QA}_f^\Lambda(\Upsilon,\nu) = |B|\int_\mathbb{R} f\left(t+\frac{\Upsilon}{2}\right)f^*\left(t-\frac{\Upsilon}{2}\right)e^{i(2At\Upsilon+B\nu t+D\Upsilon)}dt. \quad (2.8)$$

- Taking $\Lambda = (0, -1, 0, 0, 0)$, the proposed AQWD and AQAF reduces to the classical WD and AF.

Let us now examine a function $f(t) = e^{-\frac{t^2}{\sqrt{2}}} + e^{-\frac{(t-4)^2}{\sqrt{2}}}$ that sums two Gaussian beams with centers at $t = 0$ and $t = 4$, respectively. Fig. 1(a)-(b) shows the WD and its counter plot. You can see the proposed AWDQ in Fig. 1(c)-(d). Figure 2(a)-(b) shows the AF and its counter plot, whereas Figure 2(c)-(d) shows the related AAFQ. Therefore, compared to classical WD and AF, it is evident from the plots that AWDQ and AAFQ have superior detection because they can lessen the impact of cross terms on the detection process.

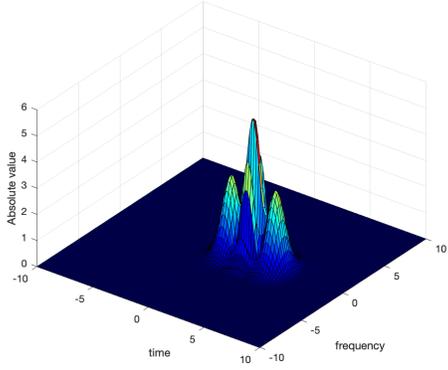
(a) WD of Gaussian beams

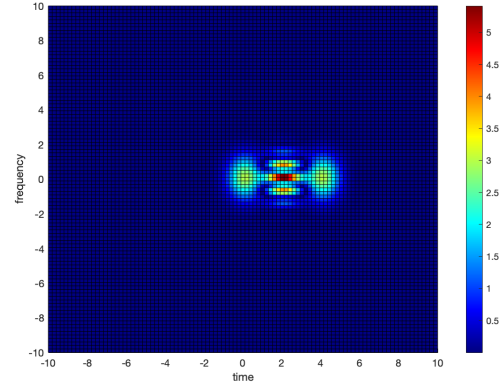
(b) Counter plot of WD of Gaussian beams .

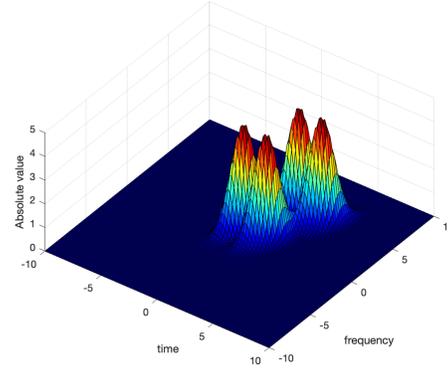
(c) AWDQ of Gaussian beams

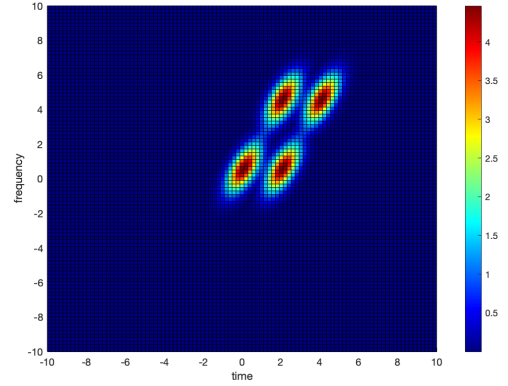
(d) Counter plot of AWD of Gaussian beams

FIGURE 1. The comparison of the absolute value of WD and AWDQ for the detection of sum of two Gaussian beams corresponding to $\Lambda = (2, -2, 0, 11)$.

Hence, it is evident from the above discussions that the proposed AWDQ and AAFQ are more adaptable than the current WD and AF classes in the QPFT domains.

## 3. Properties of the AQWD and AQAF

The primary characteristics of the suggested AQWD and AQAF will be revealed in this part. These characteristics include the Moyal formula, anti-derivative property, conjugation symmetry property, time and frequency shift properties, and time and frequency marginal features.

**(1) Shifting Properties**:

(i) Time- shifting: The AQWD and AQAF of the signal $f(t-t_0)$ have following forms:

$$\mathcal{QW}^\Lambda_{f(t-t_0)}(t,\nu) = e^{i[(C-A)(2t-t_0)t_0]} e^{i\left[2\nu + \frac{A+C}{B}t_0\right]\left(\frac{C^2-A^2}{B}t_0\right)} \times$$

$$e^{i(E-D)\left[\frac{A+C}{B}+1\right]t_0} \mathcal{QW}^\Lambda_f\left(t-t_0, \nu + \frac{(A+C)}{B}t_0\right),$$

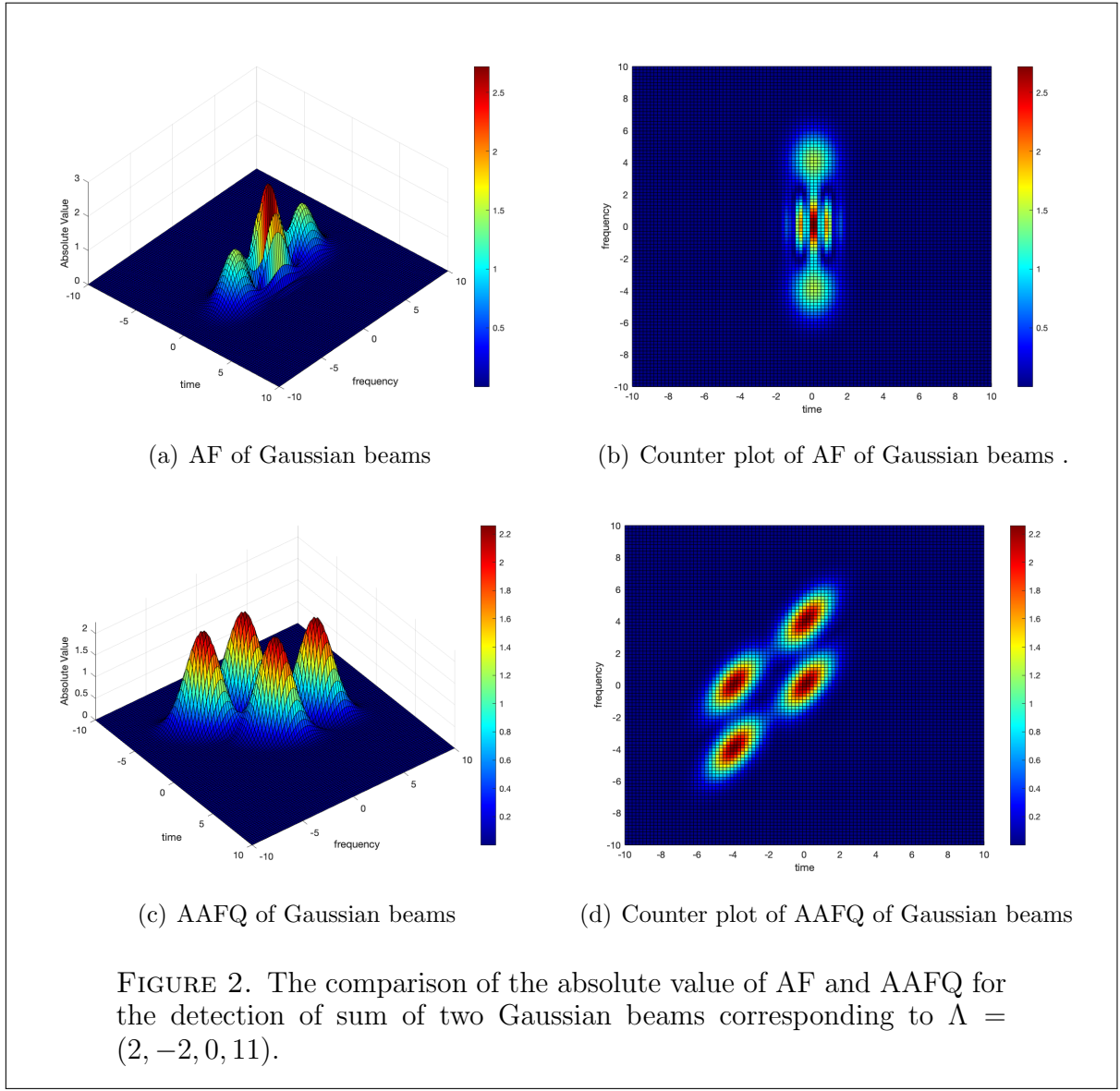

(a) AF of Gaussian beams

(b) Counter plot of AF of Gaussian beams .

(c) AAFQ of Gaussian beams

(d) Counter plot of AAFQ of Gaussian beams

FIGURE 2. The comparison of the absolute value of AF and AAFQ for the detection of sum of two Gaussian beams corresponding to $\Lambda = (2, -2, 0, 11)$.

$$\mathcal{QA}^A_{f(t-t_0)}(\Upsilon, \nu) = e^{i[(A-C)t_0 + (A+C)\Upsilon]t_0} e^{i\frac{(C-A)^2}{B}\left[\nu + \frac{C-A}{B}t_0\right]4t_0} \times$$
$$e^{i(E-D)\left[1 + \frac{(C-A)}{B}2t_0\right]} \mathcal{QA}^\Lambda_f \left(\Upsilon, \nu + \frac{(C-A)}{B}2t_0\right).$$

(ii) Frequency Shifting The AQWD and AQAF of the signal $f(t)e^{iu_0 t}$ can be given by

$$\mathcal{QW}^\Lambda_{f(t)e^{iu_0 t}}(t, \nu) = e^{i\left[(C-A)\left(\frac{u_0}{B} + 2\nu\right)\frac{u_0}{B} + (E-D)\left(\frac{u_0}{B}\right)\right]} \mathcal{QW}^\Lambda_f \left(t, \frac{u_0}{B} + \nu\right),$$
$$\mathcal{QA}^\Lambda_{f(t)e^{iu_0 t}}(\Upsilon, \nu) = e^{iu_0 \Upsilon} \mathcal{QA}^\Lambda_f(\Upsilon, \nu).$$

(iii) **Joint Time-Frequency Shifting** The AQWD and AQAF of the signal $f(t-t_0)e^{iu_0 t}$ are given by

$$\mathcal{QW}^\Lambda_{f(t-t_0)e^{iu_0 t}}(t,\nu) = e^{i[(C-A)(2t-t_0)t_0]}e^{i\left[2\nu+\frac{A+C}{B}t_0\right]\left(\frac{C^2-A^2}{B}t_0\right)}e^{i(E-D)\left[\frac{A+C}{B}+1\right]t_0} \times$$

$$e^{i\left[(C-A)\left(\frac{u_0}{B}+2\nu\right)\frac{u_0}{B}+(E-D)\left(\frac{u_0}{B}\right)\right]}\mathcal{QW}^\Lambda_f\left(t-t_0, \nu+\frac{u_0}{B}+\frac{(A+C)}{B}t_0\right),$$

$$\mathcal{QA}^\Lambda_{f(t-t_0)e^{iu_0 t}}(\Upsilon,\nu) = e^{i[(A-C)t_0+(A+C)\Upsilon]t_0}e^{i\frac{(C-A)^2}{B}\left[\nu+\frac{C-A}{B}t_0\right]4t_0} \times$$

$$e^{i(E-D)\left[1+\frac{(C-A)}{B}2t_0\right]}e^{iu_0\Upsilon}\mathcal{QA}^\Lambda_f\left(\Upsilon, \nu+\frac{(C-A)}{B}2t_0\right).$$

*Proof.* (i) The aid of (2.3) allows us to write

$$\mathcal{QW}^\Lambda_{f(t-t_0)}(t,\nu)$$

$$= |B|e^{i\left[(A-C)\nu^2+(D-E)\nu\right]}\int_\mathbb{R} f\left(t-t_0+\frac{\Upsilon}{2}\right)e^{i\left[C\left(t+\frac{\Upsilon}{2}\right)^2+\left(t+\frac{\Upsilon}{2}\right)E\right]}$$

$$\times f^*\left(t-t_0-\frac{\Upsilon}{2}\right)e^{-i\left[A\left(t-\frac{\Upsilon}{2}\right)^2+D\left(t-\frac{\Upsilon}{2}\right)\right]}e^{iB\nu\Upsilon}d\Upsilon$$

$$= |B|e^{i\left[(A-C)\nu^2+(D-E)\nu\right]}\int_\mathbb{R} f_{C,E}\left(t-t_0+\frac{\Upsilon}{2}\right)e^{i[(C-A)(2t-t_0)t_0+(E-D)t_0]}$$

$$\times f^*_{A,D}\left(t-t_0-\frac{\Upsilon}{2}\right)e^{iB\left[\nu+\frac{A+C}{B}t_0\right]\Upsilon}d\Upsilon$$

$$= e^{i[(C-A)(2t-t_0)t_0]}e^{i\left[2\nu+\frac{A+C}{B}t_0\right]\left(\frac{C^2-A^2}{B}t_0\right)}e^{i(E-D)\left[\frac{A+C}{B}+1\right]t_0}$$

$$\times|B|e^{i\left[(A-C)\left(\nu+\frac{A+C}{B}t_0\right)^2+(D-E)\left(\nu+\frac{A+C}{B}t_0\right)\right]}\int_\mathbb{R} f_{C,E}\left(t-t_0+\frac{\Upsilon}{2}\right)f^*_{A,D}\left(t-t_0-\frac{\Upsilon}{2}\right)e^{iB\left[\nu+\frac{A+C}{B}t_0\right]\Upsilon}d\Upsilon$$

$$= e^{i[(C-A)(2t-t_0)t_0]}e^{i\left[2\nu+\frac{A+C}{B}t_0\right]\left(\frac{C^2-A^2}{B}t_0\right)}e^{i(E-D)\left[\frac{A+C}{B}+1\right]t_0}\mathcal{QW}^\Lambda_f\left(t-t_0, \nu+\frac{(A+C)}{B}t_0\right)$$

Also, from (2.4), we have

$$\mathcal{QA}^\Lambda_{f(t-t_0)}(\Upsilon,\nu)$$

$$= |B|e^{i\left[(A-C)\nu^2+(D-E)\nu\right]}\int_\mathbb{R} f\left(t-t_0+\frac{\Upsilon}{2}\right)e^{i\left[C\left(t+\frac{\Upsilon}{2}\right)^2+\left(t+\frac{\Upsilon}{2}\right)E\right]}$$

$$\times f^*\left(t-t_0-\frac{\Upsilon}{2}\right)e^{-i\left[A\left(t-\frac{\Upsilon}{2}\right)^2+D\left(t-\frac{\Upsilon}{2}\right)\right]}e^{iB\nu t}dt$$

$$= |B|e^{i\left[(A-C)\nu^2+(D-E)\nu\right]}\int_\mathbb{R} f_{C,E}\left(t-t_0+\frac{\Upsilon}{2}\right)e^{i[(A-C)t_0+(E-D)+(A+C)\Upsilon]t_0}$$

$$\times f^*_{A,D}\left(t-t_0-\frac{\Upsilon}{2}\right)e^{iB\left[\nu+\frac{(C-A)}{B}2t_0\right]}dt$$

$$= e^{i[(A-C)t_0+(E-D)+(A+C)\Upsilon]t_0}e^{-i\left[(A-C)\left\{\left(\frac{(C-A)}{B}2t_0\right)^2+4\nu\frac{(C-A)}{B}t_0\right\}+(D-E)\left(\frac{(C-A)}{B}2t_0\right)\right]}$$

$$\times|B|e^{i\left[(A-C)\left(\nu+\frac{(C-A)}{B}2t_0\right)^2+(D-E)\left(\nu+\frac{(C-A)}{B}2t_0\right)\right]}\int_\mathbb{R} f_{C,E}\left(t-t_0+\frac{\Upsilon}{2}\right)f^*_{A,D}\left(t-t_0-\frac{\Upsilon}{2}\right)e^{iB\left[\nu+\frac{(C-A)}{B}2t_0\right]}d$$

$$= e^{i[(A-C)t_0+(A+C)\Upsilon]t_0}e^{i(E-D)\left[1+\frac{(C-A)}{B}2t_0\right]}e^{i\frac{(C-A)^2}{B}\left[\nu+\frac{C-A}{B}t_0\right]4t_0}\mathcal{QA}^\Lambda_f\left(\Upsilon, \nu+\frac{(C-A)}{B}2t_0\right)$$

(ii) Through basic calculations, we move forward as

$$\mathcal{QW}^{\Lambda}_{f(t)e^{iu_0 t}}(t,\nu)$$
$$= |B|e^{i\left[(A-C)\nu^2+(D-E)\nu\right]} \int_{\mathbb{R}} f_{C,E}\left(t+\frac{\Upsilon}{2}\right) e^{iu_0\left(t+\frac{\Upsilon}{2}\right)} f^*_{A,D}\left(t-\frac{\Upsilon}{2}\right) e^{-iu_0\left(t-\frac{\Upsilon}{2}\right)} e^{iB\nu\Upsilon} d\Upsilon$$
$$= |B|e^{i\left[(A-C)\nu^2+(D-E)\nu\right]} \int_{\mathbb{R}} f_{C,E}\left(t+\frac{\Upsilon}{2}\right) f^*_{A,D}\left(t-\frac{\Upsilon}{2}\right) e^{iB\left(\frac{u_0}{B}+\nu\right)\Upsilon} d\Upsilon$$
$$= e^{i\left[(C-A)\left(\frac{u_0}{B}+2\nu\right)\frac{u_0}{B}+(E-D)\left(\frac{u_0}{B}\right)\right]} \times$$
$$|B|e^{i\left[(A-C)\left(\frac{u_0}{B}+\nu\right)^2+(D-E)\left(\frac{u_0}{B}+\nu\right)\right]} \int_{\mathbb{R}} f_{C,E}\left(t+\frac{\Upsilon}{2}\right) f^*_{A,D}\left(t-\frac{\Upsilon}{2}\right) e^{iB\left(\frac{u_0}{B}+\nu\right)\Upsilon} d\Upsilon$$
$$= e^{i\left[(C-A)\left(\frac{u_0}{B}+2\nu\right)\frac{u_0}{B}+(E-D)\left(\frac{u_0}{B}\right)\right]} \mathcal{QW}^{\Lambda}_{f}\left(t,\frac{u_0}{B}+\nu\right).$$

Moreover,

$$\mathcal{QA}^{\Lambda}_{f(t)e^{iu_0 t}}(\Upsilon,\nu)$$
$$= |B|e^{i\left[(A-C)\nu^2+(D-E)\nu\right]} \int_{\mathbb{R}} f_{C,E}\left(t+\frac{\Upsilon}{2}\right) e^{iu_0\left(t+\frac{\Upsilon}{2}\right)} f^*_{A,D}\left(t-\frac{\Upsilon}{2}\right) e^{-iu_0\left(t-\frac{\Upsilon}{2}\right)} e^{iB\nu t} dt$$
$$= e^{iu_0 \Upsilon} |B|e^{i\left[(A-C)\nu^2+(D-E)\nu\right]} \int_{\mathbb{R}} f_{C,E}\left(t+\frac{\Upsilon}{2}\right) f^*_{A,D}\left(t-\frac{\Upsilon}{2}\right) e^{iB\nu t} dt$$
$$= e^{iu_0 \Upsilon} \mathcal{QA}^{\Lambda}_{f}(\Upsilon,\nu).$$

The evidence of (iii) is simple, therefore we disregard it..
This completes proof. □

**(2) Conjugation properties**:
(i) Covarriance- Conjugation

$$\left[\mathcal{QW}^{\Lambda}_{f}(t,\nu)\right]^* = -\mathcal{QW}^{\hat{\Lambda}}_{f}(t,-\nu), \qquad \left[\mathcal{QA}^{\Lambda}_{f}(\Upsilon,\nu)\right]^* = \mathcal{QA}^{\hat{\Lambda}}_{f}(-\Upsilon,-\nu),$$

where parameter $\hat{\Lambda} = (C, B, A, E, D)$, and $\tilde{\Lambda} = (C, -B, A, E, D)$.
(ii) Symmetric-Conjugation The AQWD and AQAF of the function $\check{f}(t) = f(-t)$ can be written as

$$\mathcal{QW}^{\Lambda}_{\check{f}}(t,\nu) = \mathcal{QW}^{\Lambda}_{f}(-t,-\nu), \quad \mathcal{QA}^{\Lambda}_{\check{f}}(\Upsilon,\nu) = \mathcal{QA}^{\Lambda}_{f}(-\Upsilon,-\nu).$$

*Proof.* (i) In line with (2.3), we possess

$$\left[\mathcal{QW}^{\Lambda}_{f}(t,\nu)\right]^* = |B|e^{-i\left[(A-C)\nu^2+(D-E)\nu\right]} \int_{\mathbb{R}} f^*_{C,E}\left(t+\frac{\Upsilon}{2}\right) f_{A,D}\left(t-\frac{\Upsilon}{2}\right) e^{-iB\nu\Upsilon} d\Upsilon.$$

Making change of the variable $-\Upsilon = x$, we can get the desired result as follows

$$\left[\mathcal{QW}^{\Lambda}_{f}(t,\nu)\right]^* = -|B|e^{i\left[(C-A)\nu^2+(E-D)\nu\right]} \int_{\mathbb{R}} f_{A,D}\left(t+\frac{x}{2}\right) f^*_{C,E}\left(t-\frac{x}{2}\right) e^{iB\nu x} dx$$
$$= -\mathcal{QW}^{\hat{\Lambda}}_{f}(t,\nu), \quad \text{where} \quad \hat{\Lambda} = (C, B, A, E, D).$$

Additionally, from (2.4), we obtain

$$\left[\mathcal{QA}^{\Lambda}_{f}(\Upsilon,\nu)\right]^* = |B|e^{-i\left[(A-C)\nu^2+(D-E)\nu\right]} \int_{\mathbb{R}} f^*_{C,E}\left(t+\frac{\Upsilon}{2}\right) f_{A,D}\left(t-\frac{\Upsilon}{2}\right) e^{-iB\nu t} dt$$
$$= |B|e^{i\left[(C-A)\nu^2+(E-D)\nu\right]} \int_{\mathbb{R}} f_{A,D}\left(t+\frac{-\Upsilon}{2}\right) f^*_{C,E}\left(t-\frac{-\Upsilon}{2}\right) e^{i(-B)\nu t} dt$$

$$= \mathcal{QA}_f^{\tilde{\Lambda}}(-\Upsilon, \nu), \quad \text{where} \quad \tilde{\Lambda} = (C, -B, A, E, D).$$

(ii) Additionally, the AQWD and AQAF of $\check{f}(t)$ can be shown as

$$\mathcal{QW}_{\check{f}}^{\Lambda}(t, \nu) = |B|e^{i[(A-C)\nu^2+(D-E)\nu]} \int_{\mathbb{R}} f_{C,E}\left(-t-\frac{\Upsilon}{2}\right) f_{A,D}^*\left(-t+\frac{\Upsilon}{2}\right) e^{iB\nu\Upsilon} d\Upsilon$$

$$= |B|e^{i[(A-C)\nu^2+(D-E)\nu]} \int_{\mathbb{R}} f_{C,E}\left(-t+\frac{-\Upsilon}{2}\right) f_{A,D}^*\left(-t-\frac{-\Upsilon}{2}\right) e^{i(-B)\nu(-\Upsilon)} d\Upsilon$$

$$= \mathcal{QW}_f^{\tilde{\Lambda}}(-t, \nu)$$

and

$$\mathcal{QA}_{\check{f}}^{\Lambda}(\Upsilon, \nu) = |B|e^{i[(A-C)\nu^2+(D-E)\nu]} \int_{\mathbb{R}} f_{C,E}\left(-t-\frac{\Upsilon}{2}\right) f_{A,D}^*\left(-t+\frac{\Upsilon}{2}\right) e^{iB\nu t} dt$$

$$= |B|e^{i[(A-C)\nu^2+(D-E)\nu]} \int_{\mathbb{R}} f_{C,E}\left(-t+\frac{-\Upsilon}{2}\right) f_{A,D}^*\left(-t-\frac{-\Upsilon}{2}\right) e^{i(-B)\nu(-t)} dt$$

$$= \mathcal{QA}_f^{\tilde{\Lambda}}(-\Upsilon, \nu), \quad \text{where} \quad \tilde{\Lambda} = (C, -B, A, E, D).$$

This completes proof of (ii). $\square$

**(3) QPFT marginal properties**:

(i) The following forms are found in the time and frequency marginal characteristics of the AQWD and AQAF:

$$\int_{\mathbb{R}} \mathcal{QW}_f^{\Lambda}(t, \nu) dt = 2\pi \mathbb{Q}_{\Lambda}\{f\}(\nu) \cdot \overline{\mathbb{Q}_{\Lambda'}\{f\}(\nu)}, \quad \text{where} \quad \Lambda' = (C, B, A, E, D) \quad (3.1)$$

$$\int_{\mathbb{R}} \mathcal{QA}_f^{\Lambda}(\Upsilon, \nu) d\Upsilon = 4\pi \mathbb{Q}_{\Lambda'}\{f\}(\nu) \cdot \overline{\mathbb{Q}_{\Lambda''}\{f\}(\nu)}, \quad \text{where} \quad \Lambda' = (A, B/2, C, D, E),$$

$$\text{and} \quad \Lambda'' = (C, -B/2, A, E, D). \quad (3.2)$$

(ii) The relationships shown below are valid.

$$|f(t)|^2 = \frac{1}{2\pi|B|} \int_{\mathbb{R}} e^{-i[(A-C)(\nu^2+t^2)+(D-E)(\nu+t)]} \mathcal{QW}_f^{\Lambda}(t, \nu) d\nu, \quad (3.3)$$

$$f\left(\frac{\Upsilon}{2}\right) f^*\left(-\frac{\Upsilon}{2}\right) = \frac{1}{2\pi|B|} \int_{\mathbb{R}} e^{\left[(A-C)\left(\nu^2-\left(\frac{\Upsilon}{2}\right)^2\right)+(D-E)\nu+(D+E)\frac{\Upsilon}{2}\right]} \mathcal{QA}_f^{\Lambda}(\Upsilon, \nu) d\nu. \quad (3.4)$$

*Proof.* (i) From (2.3) and (2.4), we have

$$\int_{\mathbb{R}} \mathcal{QW}_f^{\Lambda}(t, \nu) dt = |B|e^{i[(A-C)\nu^2+(D-E)\nu]} \int_{\mathbb{R}} f_{C,E}\left(t+\frac{\Upsilon}{2}\right) f_{A,D}^*\left(t-\frac{\Upsilon}{2}\right) e^{iB\nu\Upsilon} d\Upsilon dt$$

$$\int_{\mathbb{R}} \mathcal{QA}_f^{\Lambda}(\Upsilon, \nu) d\Upsilon = |B|e^{i[(A-C)\nu^2+(D-E)\nu]} \int_{\mathbb{R}} f_{C,E}\left(t+\frac{\Upsilon}{2}\right) f_{A,D}^*\left(t-\frac{\Upsilon}{2}\right) e^{iB\nu t} dt d\Upsilon.$$

Then, using $x = t + \frac{\Upsilon}{2}, y = t - \frac{\Upsilon}{2}$, we get

$$\int_{\mathbb{R}} \mathcal{QW}_f^{\Lambda}(t, \nu) dt = |B|e^{i[(A-C)\nu^2+(D-E)\nu]} \int_{\mathbb{R}} \int_{\mathbb{R}} f_{C,E}(x) f_{A,D}^*(y) e^{iB\nu(x-y)} dx dy$$

$$= |B|e^{i[(A-C)\nu^2+(D-E)\nu]} \int_{\mathbb{R}} \int_{\mathbb{R}} f(x) e^{i(Cx^2+Ex)} f^*(y) e^{-i(Ay^2+Dy)} e^{iB\nu(x-y)} dx dy$$

$$= 2\pi \left( \sqrt{\frac{B}{2\pi i}} \int_{\mathbb{R}} f(x) e^{i\left(A\nu^2 + Bx\nu + Cx^2 + D\nu + Ex\right)} dx \right) \times$$

$$\left( \sqrt{\frac{B}{-2\pi i}} \int_{\mathbb{R}} f^*(y) e^{-i\left(C\nu^2 + By\nu + Ay^2 + E\nu + Dy\right)} dy \right)$$

$$= 2\pi \left( \sqrt{\frac{B}{2\pi i}} \int_{\mathbb{R}} f(x) e^{i\left(A\nu^2 + Bx\nu + Cx^2 + D\nu + Ex\right)} dx \right) \times$$

$$\left( \sqrt{\frac{B}{2\pi i}} \int_{\mathbb{R}} f(y) e^{i\left(C\nu^2 + By\nu + Ay^2 + E\nu + Dy\right)} dy \right)^*$$

$$= 2\pi \mathbb{Q}_\Lambda\{f\}(\nu) \cdot \overline{\mathbb{Q}_{\Lambda'}\{f\}(\nu)}, \quad \text{where} \quad \Lambda' = (C, B, A, E, D).$$

In addition

$$\int_{\mathbb{R}} \mathcal{Q}\mathcal{A}_f^\Lambda(\Upsilon, \nu) d\Upsilon = |B| e^{i\left[(A-C)\nu^2 + (D-E)\nu\right]} \int_{\mathbb{R}} \int_{\mathbb{R}} f_{C,E}(x) f_{A,D}^*(y) e^{iB\nu \frac{(x+y)}{2}} dx dy$$

$$= |B| e^{i\left[(A-C)\nu^2 + (D-E)\nu\right]} \int_{\mathbb{R}} \int_{\mathbb{R}} f(x) e^{i(Cx^2 + Ex)} f^*(y) e^{-i(Ay^2 + Dy)} e^{iB\nu \frac{(x+y)}{2}} dx dy$$

$$= 4\pi \left( \sqrt{\frac{B/2}{2\pi i}} \int_{\mathbb{R}} f(x) e^{i\left(A\nu^2 + \frac{B}{2} x\nu + Cx^2 + D\nu + Ex\right)} dx \right) \times$$

$$\left( \sqrt{\frac{-B/2}{-2\pi i}} \int_{\mathbb{R}} f^*(y) e^{-i\left(C\nu^2 + \frac{-B}{2} y\nu + Ay^2 + E\nu + Dy\right)} dy \right)$$

$$= 4\pi \left( \sqrt{\frac{B/2}{2\pi i}} \int_{\mathbb{R}} f(x) e^{i\left(A\nu^2 + \left(\frac{B}{2}\right) x\nu + Cx^2 + D\nu + Ex\right)} dx \right) \times$$

$$\left( \sqrt{\frac{-B/2}{2\pi i}} \int_{\mathbb{R}} f(y) e^{i\left(C\nu^2 + \left(\frac{-B}{2}\right) y\nu + Ay^2 + E\nu + Dy\right)} dy \right)^*$$

$$= 4\pi \mathbb{Q}_{\Lambda'}\{f\}(\nu) \cdot \overline{\mathbb{Q}_{\Lambda''}\{f\}(\nu)}, \quad \text{where} \quad \Lambda' = (A, B/2, C, D, E),$$
$$\text{and} \quad \Lambda'' = (C, -B/2, A, E, D).$$

Hence, (3.1) and (3.2) are proved.

(ii) From (2.3) and (2.4), we get

$$\frac{1}{|B|} e^{-i\left[(A-C)\nu^2 + (D-E)\nu\right]} \mathcal{Q}\mathcal{W}_f^\Lambda(t, \nu) = \int_{\mathbb{R}} f_{C,E}\left(t + \frac{\Upsilon}{2}\right) f_{A,D}^*\left(t - \frac{\Upsilon}{2}\right) e^{iB\nu\Upsilon} d\Upsilon,$$

$$\frac{1}{|B|} e^{-i\left[(A-C)\nu^2 + (D-E)\nu\right]} \mathcal{Q}\mathcal{A}_f^\Lambda(\Upsilon, \nu) = \int_{\mathbb{R}} f_{C,E}\left(t + \frac{\Upsilon}{2}\right) f_{A,D}^*\left(t - \frac{\Upsilon}{2}\right) e^{iB\nu t} dt,$$

Therefore using inverse WD and AF, we have

$$f_{C,E}\left(t + \frac{\Upsilon}{2}\right) f_{A,D}^*\left(t - \frac{\Upsilon}{2}\right) = \frac{1}{2\pi |B|} \int_{\mathbb{R}} e^{-i\left[(A-C)\nu^2 + (D-E)\nu\right]} \mathcal{Q}\mathcal{W}_f^\Lambda(t, \nu) e^{-iB\nu\Upsilon} d\nu, \quad (3.5)$$

$$f_{C,E}\left(t + \frac{\Upsilon}{2}\right) f_{A,D}^*\left(t - \frac{\Upsilon}{2}\right) = \frac{1}{2\pi |B|} \int_{\mathbb{R}} e^{-i\left[(A-C)\nu^2 + (D-E)\nu\right]} \mathcal{Q}\mathcal{A}_f^\Lambda(\Upsilon, \nu) e^{-iB\nu t} d\nu. \quad (3.6)$$

By changing $t = 0$ in (3.6) and $\Upsilon = 0$ in (3.5), we get

$$f_{C,E}(t) f_{A,D}^*(t) = \frac{1}{2\pi|B|} \int_{\mathbb{R}} e^{-i[(A-C)\nu^2 + (D-E)\nu]} \mathcal{QW}_f^\Lambda(t,\nu) d\nu,$$

$$f_{C,E}\left(\frac{\Upsilon}{2}\right) f_{A,D}^*\left(-\frac{\Upsilon}{2}\right) = \frac{1}{2\pi|B|} \int_{\mathbb{R}} e^{-i[(A-C)\nu^2 + (D-E)\nu]} \mathcal{QA}_f^\Lambda(\Upsilon,\nu) d\nu.$$

The aforementioned equations can be rewritten as

$$|f(t)|^2 = \frac{1}{2\pi|B|} \int_{\mathbb{R}} e^{-i[(A-C)(\nu^2+t^2)+(D-E)(\nu+t)]} \mathcal{QW}_f^\Lambda(t,\nu) d\nu,$$

$$f\left(\frac{\Upsilon}{2}\right) f^*\left(-\frac{\Upsilon}{2}\right) = \frac{1}{2\pi|B|} \int_{\mathbb{R}} e^{\left[(A-C)\left(\nu^2-\left(\frac{\Upsilon}{2}\right)^2\right)+(D-E)\nu+(D+E)\frac{\Upsilon}{2}\right]} \mathcal{QA}_f^\Lambda(\Upsilon,\nu) d\nu.$$

This, proves (3.3) and (3.4).
This completes proof. □

**(4) Moyal's formula**:

For the AQWD and AQAF, the Moyal's formula is provided by

$$\int_{\mathbb{R}} \int_{\mathbb{R}} \mathcal{QW}_f^\Lambda(t,\nu) [\mathcal{QW}_g^\Lambda(t,\nu)]^* dt d\nu = 2\pi B |\langle f, g \rangle|^2, \tag{3.7}$$

$$\int_{\mathbb{R}} \int_{\mathbb{R}} \mathcal{QA}_f^\Lambda(\Upsilon,\nu) [\mathcal{QA}_g^\Lambda(\Upsilon,\nu)]^* d\Upsilon d\nu = 2\pi B |\langle f, g \rangle|^2, \tag{3.8}$$

where $f, g \in L^2(\mathbb{R})$.

*Proof.* To prove (3.7), we begin with

$$\int_{\mathbb{R}} \int_{\mathbb{R}} \mathcal{QW}_f^\Lambda(t,\nu) [\mathcal{QW}_g^\Lambda(t,\nu)]^* dt d\nu$$

$$= B^2 \int_{\mathbb{R}^4} f_{C,E}\left(t + \frac{\Upsilon}{2}\right) f_{A,D}^*\left(t - \frac{\Upsilon}{2}\right) g_{C,E}^*\left(t + \frac{\varepsilon}{2}\right) g_{A,D}\left(t - \frac{\varepsilon}{2}\right) e^{-iB\nu(\varepsilon-\Upsilon)} d\Upsilon d\varepsilon dt d\nu$$

$$= 2\pi|B| \int_{\mathbb{R}^3} f_{C,E}\left(t + \frac{\Upsilon}{2}\right) f_{A,D}^*\left(t - \frac{\Upsilon}{2}\right) g_{C,E}^*\left(t + \frac{\varepsilon}{2}\right) g_{A,D}\left(t - \frac{\varepsilon}{2}\right) \left(\frac{B}{2\pi} \int_{\mathbb{R}} e^{iB\nu(\Upsilon-\varepsilon)} d\nu\right) d\Upsilon d\varepsilon dt$$

$$= 2\pi B \int_{\mathbb{R}^3} f_{C,E}\left(t + \frac{\Upsilon}{2}\right) f_{A,D}^*\left(t - \frac{\Upsilon}{2}\right) g_{C,E}^*\left(t + \frac{\varepsilon}{2}\right) g_{A,D}\left(t - \frac{\varepsilon}{2}\right) \delta(\Upsilon - \varepsilon) d\varepsilon d\Upsilon dt$$

$$= 2\pi B \int_{\mathbb{R}^2} f_{C,E}\left(t + \frac{\Upsilon}{2}\right) f_{A,D}^*\left(t - \frac{\Upsilon}{2}\right) g_{C,E}^*\left(t + \frac{\Upsilon}{2}\right) g_{A,D}\left(t - \frac{\Upsilon}{2}\right) d\Upsilon dt.$$

Making the substitution $u = t + \dfrac{\Upsilon}{2}$ and $v = t - \dfrac{\Upsilon}{2}$ above equation yields

$$\int_{\mathbb{R}} \int_{\mathbb{R}} \mathcal{QW}_f^{\Lambda}(\Upsilon, \nu)[\mathcal{QW}_g^{\Lambda}(\Upsilon, \nu)]^* dt d\nu$$

$$= 2\pi B \int_{\mathbb{R}^2} f_{C,E}(u) f_{A,D}^*(v) g_{C,E}^*(u) g_{A,D}(v) \, du dv$$

$$= 2\pi B \int_{\mathbb{R}} f_{C,E}(u) g_{C,E}^*(u) \, du \cdot \int_{\mathbb{R}} f_{A,D}^*(v) g_{A,D}(v) \, dv$$

$$= 2\pi B \left( \int_{\mathbb{R}} f(u) g^*(u) \, du \right) \cdot \left( \int_{\mathbb{R}} f(v) g^*(v) \, dv \right)^*$$

$$= 2\pi B |\langle f, g \rangle|^2,$$

It is the intended outcome.

Given how closely the proof of (3.8) resembles that of (3.7), it will be excluded.

Hence completes the proof.

□

**(5) Anti-derivative property**:

The following formulas can be used to reconstruct $f \in L^2(\mathbb{R})$ from the proposed AQWD and AQAF:

$$f(t) = \frac{e^{i[C\nu^2 + E\nu]}}{2\pi f^*(0)|B|\sqrt{iB}} \int_{\mathbb{R}} \mathcal{K}_\Lambda^*(\nu, t) \, \mathcal{QW}_f^\Lambda\left(\frac{t}{2}, \nu\right) d\nu, \tag{3.9}$$

$$f(t) = \frac{e^{i[C\nu^2 + E\nu]}}{2\pi f^*(0)|B|\sqrt{iB}} \int_{\mathbb{R}} \mathcal{K}_\Lambda^*\left(\frac{\nu}{2}, t\right) \mathcal{QA}_f^\Lambda(2t, \nu) \, d\nu, \tag{3.10}$$

provided $f^*(0) \neq 0$.

*Proof.* First, we demonstrate (3.9). On taking $\dfrac{\Upsilon}{2}$ with $t$ in (3.5), we get

$$f_{C,E}(2t) f_{A,D}^*(0) = \frac{1}{2\pi |B|} \int_{\mathbb{R}} e^{-i[(A-C)\nu^2 + (D-E)\nu]} \mathcal{QW}_f^\Lambda(t, \nu) \, e^{-iB\nu(2t)} d\nu.$$

On further simplification, it yields

$$f(2t) f^*(0) = \frac{1}{2\pi |B|} e^{-i[C(2t)^2 + E(2t)]} \int_{\mathbb{R}} e^{-i[(A-C)\nu^2 + (D-E)\nu]} \mathcal{QW}_f^\Lambda(t, \nu) \, e^{-iB(2t)\nu} d\nu$$

$$= \frac{1}{2\pi |B|} e^{i[C\nu^2 + E\nu]} \int_{\mathbb{R}} e^{-i[A\nu^2 + B(2t)\nu + C(2t)^2 + D\nu + E(2t)]} \mathcal{QW}_f^\Lambda(t, \nu) \, d\nu$$

$$= \frac{e^{i[C\nu^2 + E\nu]}}{2\pi |B|\sqrt{iB}} \int_{\mathbb{R}} \left( \sqrt{\frac{B}{i}} e^{i[A\nu^2 + B(2t)\nu + C(2t)^2 + D\nu + E(2t)]} \right)^* \mathcal{QW}_f^\Lambda(t, \nu) \, d\nu.$$

The above formula becomes

$$f(t) = \frac{e^{i[C\nu^2 + E\nu]}}{2\pi f^*(0)|B|\sqrt{iB}} \int_{\mathbb{R}} \mathcal{K}_\Lambda^*(\nu, t) \, \mathcal{QW}_f^\Lambda\left(\frac{t}{2}, \nu\right) d\nu.$$

We obtain (3.10) in a manner similar to the proof of (3.9).

Thus, the proof is finished.

□

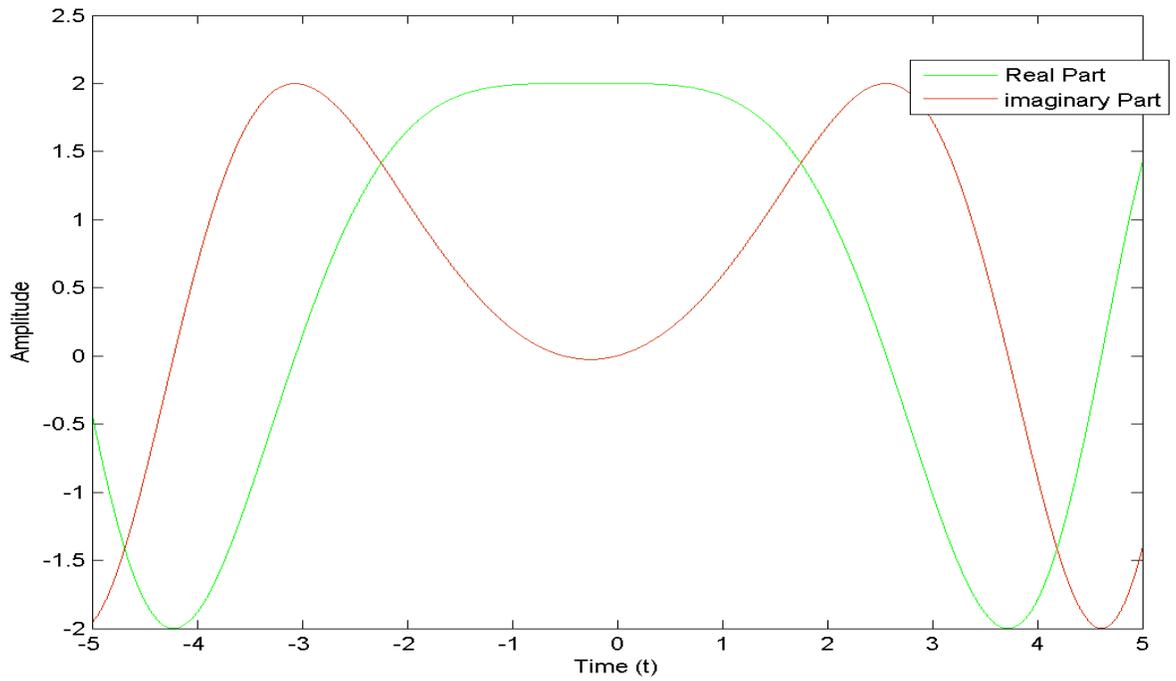

FIGURE 3. Real and Imaginary parts of mono-component signal $r(t) = e^{i(0.1t+0.2t^2)}$

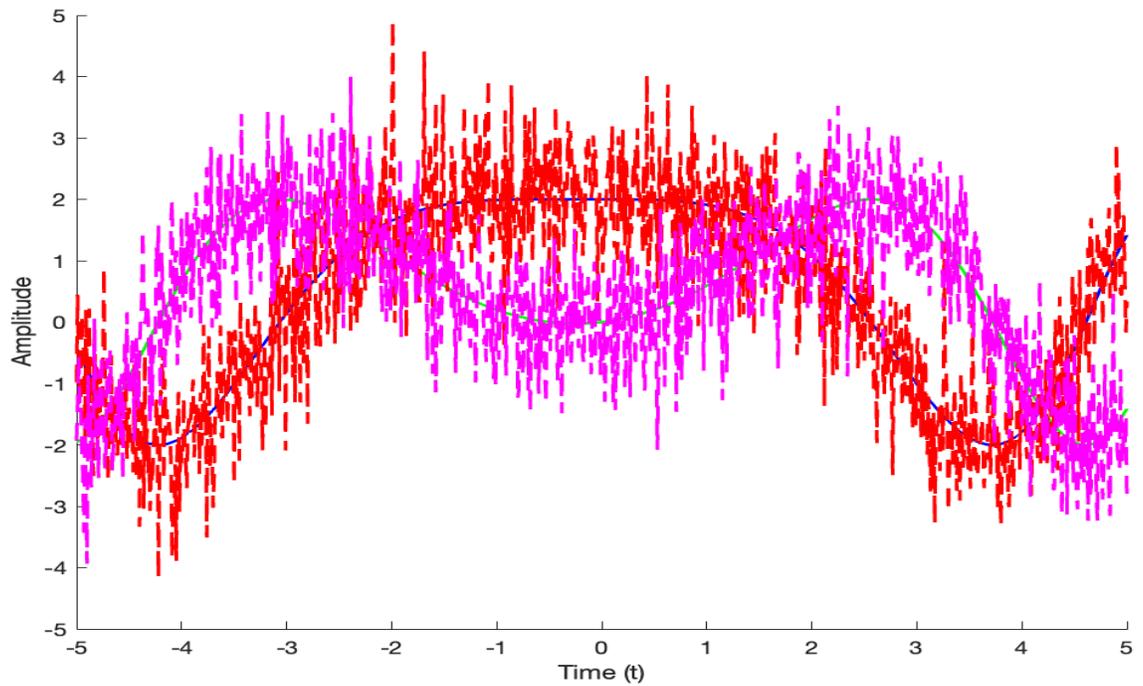

FIGURE 4. Real and Imaginary parts of mono-component signal $r(t) = e^{i(0.1t+0.2t^2)}$, at 5dB SNR.

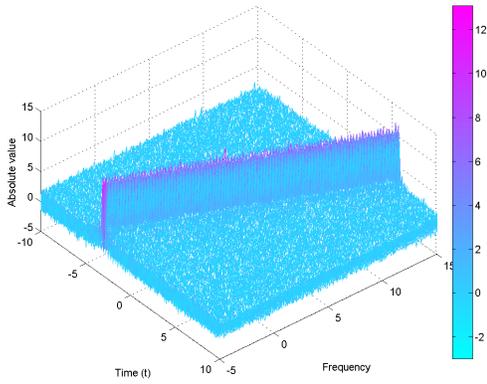
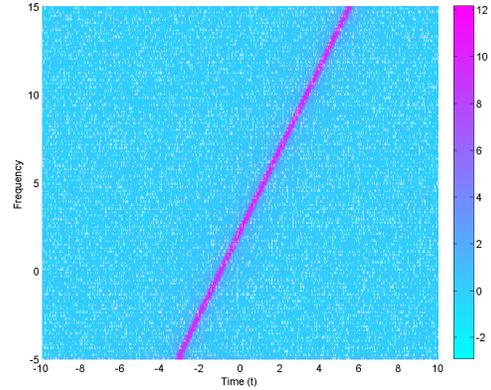

(a) AQWD of $r(t)$ with SNR = 5dB.

(b) Contour plot of AQWD of $r(t)$ with SNR = 5dB.

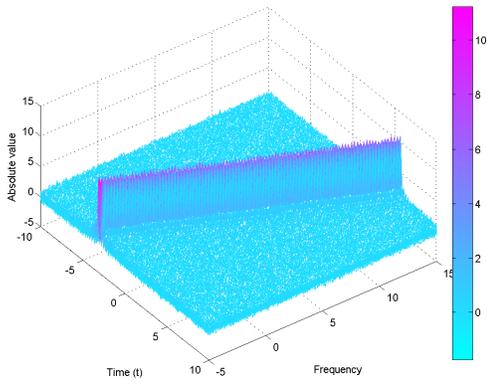
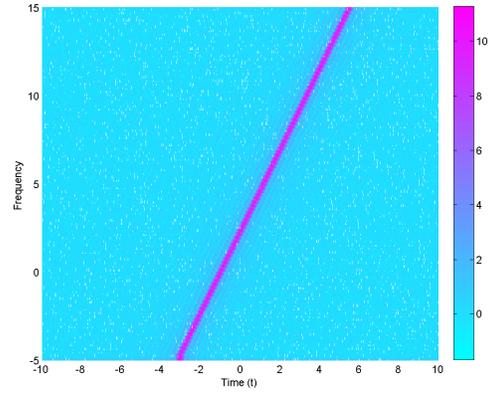

(c) AQWDof $r(t)$ with SNR = 10dB.

(d) Contour plot of AQWD of $r(t)$ with SNR = 10dB.

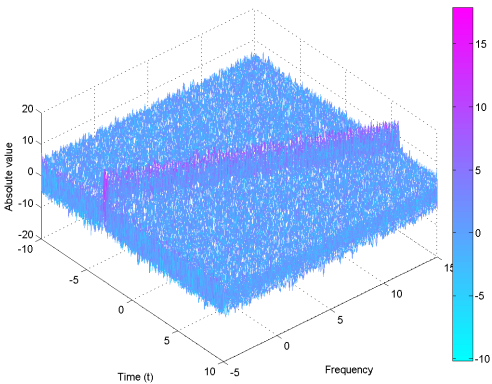
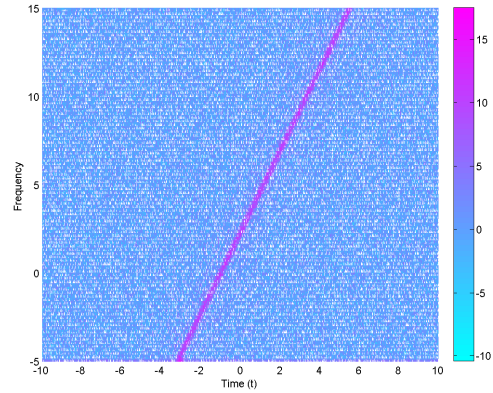

(e) AQWD of $r(t)$ with SNR = -5dB.

(f) Contour plot of AQWD of $r(t)$ with SNR = -5dB.

FIGURE 5. The detection and parameters estimation for $r(t) = e^{i(0.1t+0.2t^2)}$ with $\Lambda = (0, -1, 0, 2, 2)$ and noise by the AQWD.

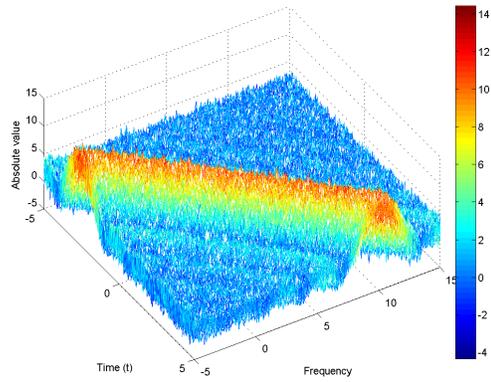

(a) WD of $r(t)$.

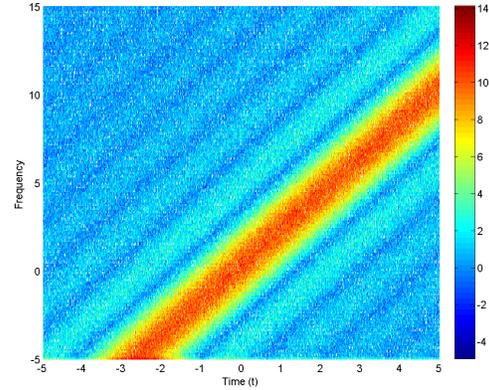

(b) Contour plot of WD of $r(t)$.

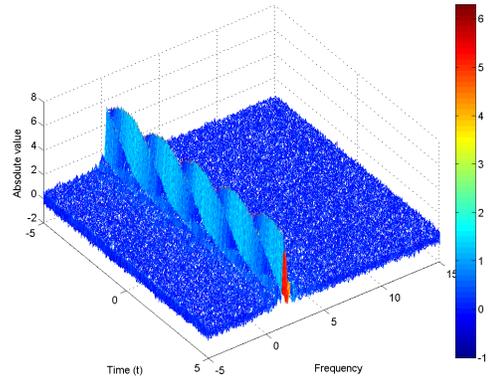

(c) QWD of $r(t)$ with $\Lambda = (0, -2, 1, 2, 1)$.

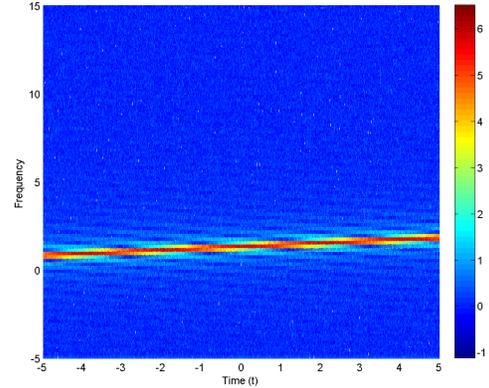

(d) Contour plot of QWD of $r(t)$ with $\Lambda = (0, -2, 1, 2, 1)$.

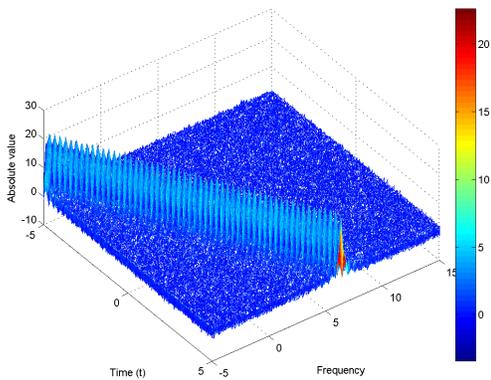

(e) AQWD of $r(t)$ with $\Lambda = (1, -2, 1, 2, 1)$.

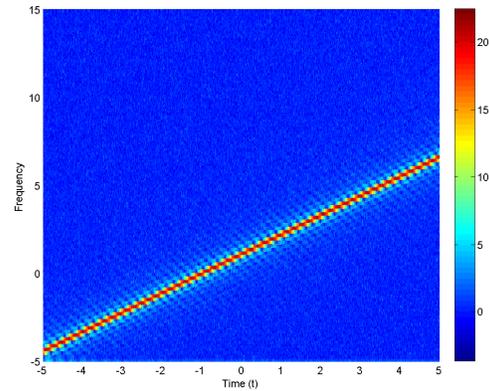

(f) Contour plot of AQWD of $r(t)$ with $\Lambda = (1, -2, 1, 2, 1)$.

FIGURE 6. The comparison of the absolute value of WD, QWD and AQWD for the detection of mono-component LFM signal $r(t) = e^{i(0.1t + 0.2t^2)}$ corresponding to specific $\Lambda$ options with $A_0 = 1$, $T = 10$ and SNR $=10$dB.

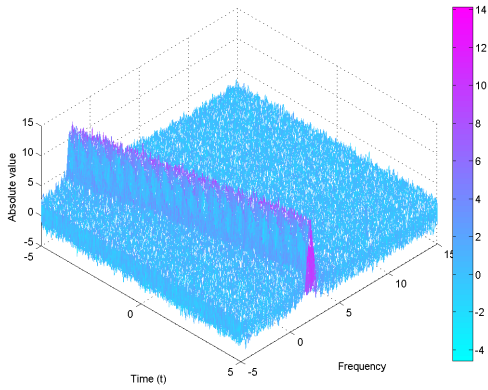

(a) AQAF of $r(t)$ with SNR = 5dB.

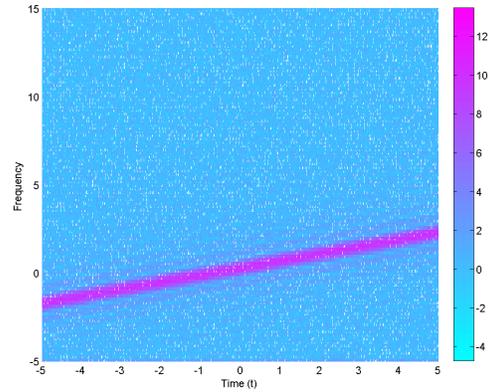

(b) Contour plot of AQAF of $r(t)$ with SNR = 5dB.

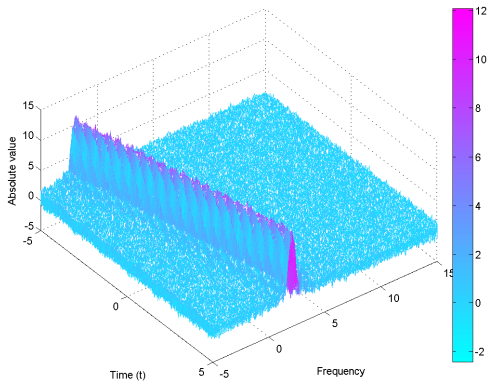

(c) AQAF of $r(t)$ with SNR = 10dB.

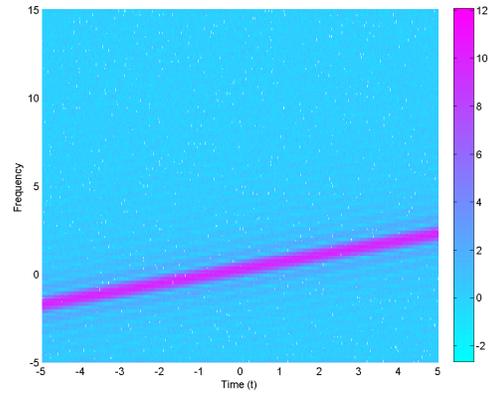

(d) Contour plot of AQAF of $r(t)$ with SNR = 10dB.

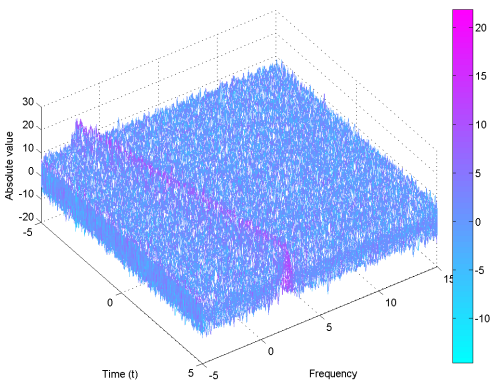

(e) AQAF of $r(t)$ with SNR = -5dB.

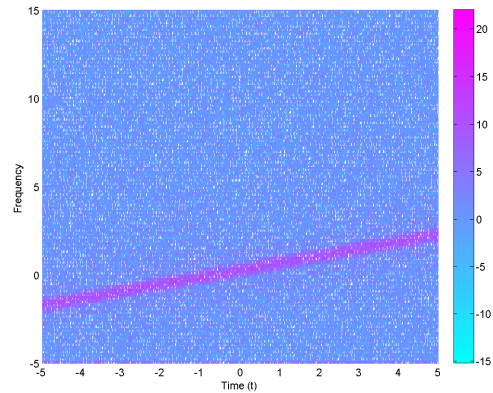

(f) Contour plot of AQAF of $r(t)$ with SNR = -5dB.

FIGURE 7. The detection and parameters estimation for $r(t) = e^{i(0.1t+0.2t^2)}$ with $\Lambda = (0, -1, 0, 2, 2)$ and noise by the AQAF.

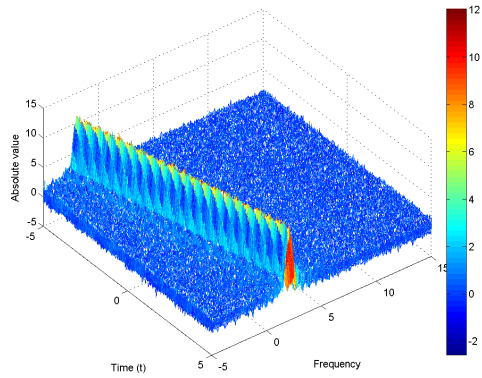

(a) AF of $r(t)$.

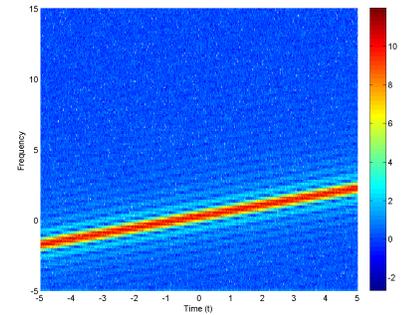

(b) Contour plot of AF of $r(t)$.

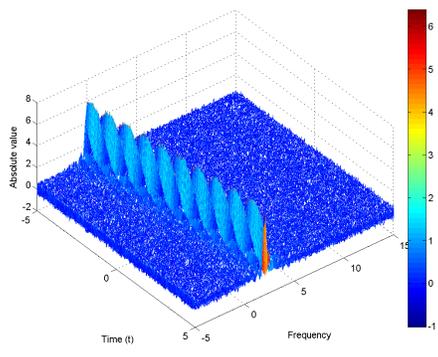

(c) QAF of $r(t)$ with $\Lambda = (0, -2, 1, 2, 1)$.

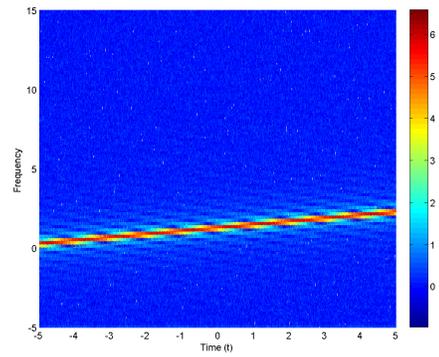

(d) Contour plot of QAF of $r(t)$ with $\Lambda = (0, -2, 1, 2, 1)$.

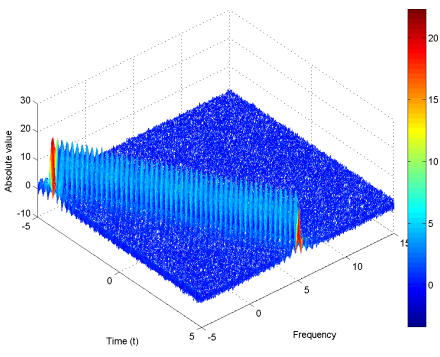

(e) AQAF of $r(t)$ with $\Lambda = (1, -2, 1, 2, 1)$.

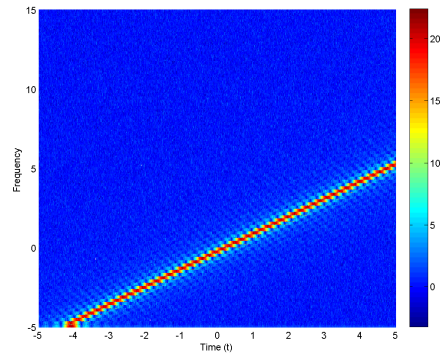

(f) Contour plot of AQAF of $r(t)$ with $\Lambda = (1, -2, 1, 2, 1)$.

FIGURE 8. The comparison of the absolute value of AF, QAF and AQAF for the detection of one-component LFM signal $r(t) = e^{i(0.1t + 0.2t^2)}$ corresponding to specific values of $\Lambda$ with $A_0 = 1$, $T = 10$ and SNR =10dB.

## 4. Applications

It is common to come across LFM signals, which are typical non-stationary signals, in applications like communications and sonar [18]. This section will delve into a detailed investigation of the applications of AQWD and AQAF in the detection of both single and multi-component LFM signals.

4.1. **Single(mono)-component LFM signal.** Examine the LFM signal $f(t)$, which is a single component, as follows:

$$f(t) = A_0 e^{i(\nu_0 t + \xi_0 t^2)}, \quad -\frac{T}{2} \leq t \leq \frac{T}{2},$$

where the amplitude $A_0$, initial frequency $\nu_0$, and frequency rate $\xi_0$. The AQWD of $f(t)$ can be find as

$$\mathcal{QW}_f^\Lambda(t,\nu)$$
$$= |B| e^{i\left[(A-C)\nu^2 + (D-E)\nu\right]} \int_{-\frac{T}{2}}^{\frac{T}{2}} A_0 e^{i\left[\nu_0\left(t+\frac{\Upsilon}{2}\right) + \xi_0\left(t+\frac{\Upsilon}{2}\right)^2\right]} e^{-i\left[\nu_0\left(t-\frac{\Upsilon}{2}\right) + \xi_0\left(t-\frac{\Upsilon}{2}\right)^2\right]}$$
$$\times A_0^* e^{i\left[C\left(t+\frac{\Upsilon}{2}\right)^2 + E\left(t+\frac{\Upsilon}{2}\right)\right]} e^{-i\left[A\left(t-\frac{\Upsilon}{2}\right)^2 + D\left(t-\frac{\Upsilon}{2}\right)\right]} e^{iB\nu\Upsilon} d\Upsilon$$
$$= |A_0|^2 |B| e^{i\left[(A-C)(\nu^2 - t^2) + (D-E)(\nu - t)\right]} \int_{-\frac{T}{2}}^{\frac{T}{2}} e^{i\frac{(C-A)}{4}\Upsilon^2} e^{i\left[(2\xi_0 + C + A)t + B\nu + \nu_0 + \frac{D+E}{2}\right]\Upsilon} d\Upsilon.$$

Hence

$$\mathcal{QW}_f^\Lambda(t,\nu)$$
$$= \begin{cases} |A_0|^2 T |B| e^{i[(D-E)(\nu-t)]} \mathrm{sinc}\left\{\frac{T}{2}\left[(2\xi_0 + 2A)t + B\nu + \nu_0 + \frac{D+E}{2}\right]\right\}, & A = C \\ |A_0|^2 |B| e^{i\left[(A-C)(\nu^2 - t^2) + (D-E)(\nu - t)\right]} \int_{-\frac{T}{2}}^{\frac{T}{2}} e^{i\frac{(C-A)}{4}\Upsilon^2} e^{i\left[(2\xi_0 + C + A)t + B\nu + \nu_0 + \frac{D+E}{2}\right]\Upsilon} d\Upsilon, & A \neq C. \end{cases}$$
(4.1)

or

$$\mathcal{QW}_f^A(t,\nu)$$
$$= \begin{cases} |A_0|^2 T |B| \mathrm{sinc}\left\{\frac{T}{2}\left[(2\xi_0 + 2A)t + B\nu + \nu_0 + D\right]\right\}, & A = C, D = E \\ |A_0|^2 |B| e^{i\left[(A-C)(\nu^2 - t^2) + (D-E)(\nu - t)\right]} \int_{-\frac{T}{2}}^{\frac{T}{2}} e^{i\frac{(C-A)}{4}\Upsilon^2} e^{i\left[(2\xi_0 + C + A)t + B\nu + \nu_0 + \frac{D+E}{2}\right]\Upsilon} d\Upsilon, & A \neq C, D \neq E. \end{cases}$$
(4.2)

The AQWD of a single-component LFM signal $f(t)$ can be inferred from (4.1) to generate impulses at a straight line $\left[(2\xi_0 + 2A)t + B\nu + \nu_0 + \frac{D+E}{2}\right] = 0$ in the $(t,\nu)$-plane, where $A = C$. Thus, we can conclude that by appropriately selecting the parameters, the AQAF can be used to detect a single-component LFM signal. Additionally, the QWD of a single-component LFM signal $f(t)$ can be obtained using (1.8) by

$$\mathrm{QWD}_f^A(\Upsilon, \nu) = \begin{cases} |A_0|^2 \sqrt{\frac{B}{2\pi i}} e^{i(C\nu^2 + E\nu)} T \mathrm{sinc}\left\{\frac{T}{2}\left[2\xi_0 t + B\nu + \nu_0 + D\right]\right\}, & A = 0 \\ |A_0|^2 \sqrt{\frac{B}{2\pi i}} e^{i(C\nu^2 + E\nu)} \int_{-\frac{T}{2}}^{\frac{T}{2}} e^{i\left(A\Upsilon^2 + i\Upsilon(B\nu + \nu_0 + 2\xi_0 t + D)\right)} d\Upsilon, & A \neq 0. \end{cases}$$
(4.3)

However, in the $(t,\nu)$-plane, the QWD is unable to produce impulses at a straight line when $A \neq 0$, but the proposed AQWD is capable of doing so. Thus, in order to detect a single-component LFM, it is more convenient to employ the AQWD rather than the QWD.

For example, the detection and estimate for the single-component LFM signal $r(t) = e^{i(0.1t+0.2t^2)}$ with SNR = 5dB, SNR = 10dB, and SNR = -5 by the AQWD for the parameters $\Lambda = (0, -1, 0, 2, 2)$ are shown in Fig. 5. Furthermore, it can be observed from Fig. 5 that the contour images of AQWD can be used to detect LFM signal. Additionally, Fig. 6 compares the parameters estimate and detection for $r(t) = e^{i(0.1t+0.2t^2)}$ ($|t| \leq 10$) with SNR = 10dB using WD, QWD, and AQWD. According to their contour images, AQWD appears to be a more successful method for detecting a single-component LFM than WD and QWD.

Similarly, the AQAF of the mono-component signal $f(t)$ can be found as

$$\mathcal{QA}_f^\Lambda(t, \nu)$$
$$= |B|e^{i[(A-C)\nu^2 + (D-E)\nu]} \int_{-\frac{T}{2}}^{\frac{T}{2}} A_0 e^{i\left[\nu_0\left(t+\frac{\Upsilon}{2}\right) + \xi_0\left(t+\frac{\Upsilon}{2}\right)^2\right]} e^{-i\left[\nu_0\left(t-\frac{\Upsilon}{2}\right) + \xi_0\left(t-\frac{\Upsilon}{2}\right)^2\right]}$$
$$\times A_0^* e^{i\left[C\left(t+\frac{\Upsilon}{2}\right)^2 + E\left(t+\frac{\Upsilon}{2}\right)\right]} e^{-i\left[A\left(t-\frac{\Upsilon}{2}\right)^2 + D\left(t-\frac{\Upsilon}{2}\right)\right]} e^{iB\nu t} dt$$
$$= |A_0|^2 |B| e^{i\left[(A-C)(\nu^2 - \Upsilon^2/4) + (D-E)(\nu - \Upsilon/2)\right]} e^{i\nu_0 \Upsilon} \int_{-\frac{T}{2}}^{\frac{T}{2}} e^{i(C-A)t^2} e^{i[(2\xi_0 + C+A)\Upsilon + B\nu + E - D]t} dt.$$

Hence,

$$\mathcal{QA}_f^\Lambda(\Upsilon, \nu)$$
$$= \begin{cases} |A_0|^2 |B| e^{i[(D-E)(\nu-\Upsilon/2)]} e^{i\nu_0\Upsilon} T \text{sinc}\left\{\frac{T}{2}\left[(2\xi_0 + 2A)\Upsilon + B\nu + E - D\right]\right\}, & A = C \\ |A_0|^2 |B| e^{i[(A-C)(\nu^2-\Upsilon^2/4)+(D-E)(\nu-\Upsilon/2)]} e^{i\nu_0\Upsilon} \int_{-\frac{T}{2}}^{\frac{T}{2}} e^{i(C-A)t^2} e^{i[(2\xi_0+C+A)\Upsilon+B\nu+E-D]t} dt, A \neq C. \end{cases}$$

Or

$$\mathcal{QA}_f^\Lambda(\Upsilon, \nu)$$
$$= \begin{cases} |A_0|^2 |B| e^{i\nu_0\Upsilon} T \text{sinc}\left\{\frac{T}{2}\left[(2\xi_0 + 2A)\Upsilon + B\nu\right]\right\}, & A = C, D = E \\ |A_0|^2 |B| e^{i[(A-C)(\nu^2-\Upsilon^2/4)+(D-E)(\nu-\Upsilon/2)]} e^{i\nu_0\Upsilon} \int_{-\frac{T}{2}}^{\frac{T}{2}} e^{i(C-A)t^2} e^{i[(2\xi_0+C+A)\Upsilon+B\nu+E-D]t} dt, \\ & A \neq C, D \neq E. \end{cases}$$

Consequently, the AQAF of a single-component LFM signal $f(t)$ produces impulses at a straight line $[(2\xi_0 + 2A)\Upsilon + B\nu + E - D] = 0$ when $A = C$ in the $(\Upsilon, \nu)$ plane. Thus, we can conclude that by appropriately selecting the parameters, the AQAF can be used to detect a single-component LFM signal. Additionally, we extract the QAF of the single-component LFM signal $f(t)$ using (1.9) as follows:

$$\text{QAF}_f^\Lambda(\Upsilon, \nu) = \begin{cases} |A_0|^2 \sqrt{\frac{B}{2\pi i}} e^{i(C\nu^2 + E\nu)} e^{i\nu_0\Upsilon} T \text{sinc}\left\{\frac{T}{2}\left[2\xi_0 t + B\nu + D\right]\right\}, & A = 0 \\ |A_0|^2 \sqrt{\frac{B}{2\pi i}} e^{i(C\nu^2 + E\nu)} e^{i\nu_0\Upsilon} \int_{-\frac{T}{2}}^{\frac{T}{2}} e^{i\left(At^2 + i(B\nu + 2\xi_0 t + D)t\right)} dt, & A \neq 0. \end{cases}$$

However, in the $(\Upsilon, \nu)$-plane, the QAF is unable to produce impulses at a straight line when $A \neq 0$, but the proposed AQAF is capable of doing so. Thus, we conclude that the new AQAF is more flexible than QAF in the detection of a single-component LFM signal.

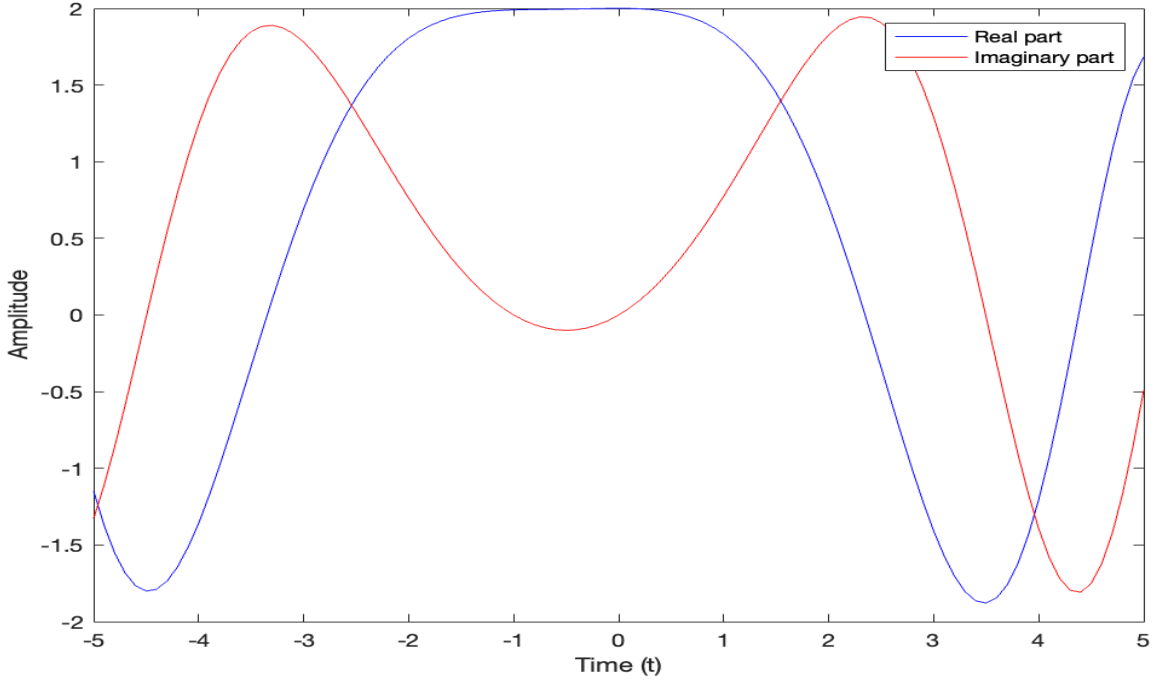

FIGURE 9. Real and Imaginary parts of bi-component signal $u(t) = e^{i(0.1t+0.2t^2)} + e^{i(0.3t+0.2t^2)}$

For example, the detection and estimate for the single-component LFM signal $r(t) = e^{i(0.1t+0.2t^2)}$ with SNR = 5dB, SNR = 10dB, and SNR = -5 by the AQAF for the parameters $\Lambda = (0, -1, 0, 2, 2)$ are shown in Fig. 7. Furthermore, it can be observed from Fig. 7 that the contour images of AQAF can be used to detect LFM signal. Additionally, Fig. 8 compares the parameters estimate and detection for $r(t) = e^{i(0.1t+0.2t^2)}$ ($|t| \leq 10$) with SNR = 10dB using AF, QAF, and AQAF. Therefore, using AQAF will be significantly more effective than using AF and QAF in detecting single-component LFM signals.

4.2. **Multi-component LFM signal.** The general form of multi-component LFM signal is given by

$$f(t) = \sum_{k=1}^{n} f_k(t), \quad \frac{T}{2} \leq t \leq \frac{T}{2},$$

where $f_k(t) = A_k e^{i(\nu_k t + \xi_k t^2)}$, $k \in \mathbb{N}$.

The AQWD of $f(t)$ can be expressed as

$$\mathcal{QW}_f^\Lambda(t,\nu) = \sum_{k=1}^{n} \mathcal{QW}_{f_k}^\Lambda(t,\nu) + \sum_{k_1 \neq k_2 = 1}^{n} \mathcal{QW}_{f_{k_1}, f_{k_2}}^\Lambda(t,\nu).$$

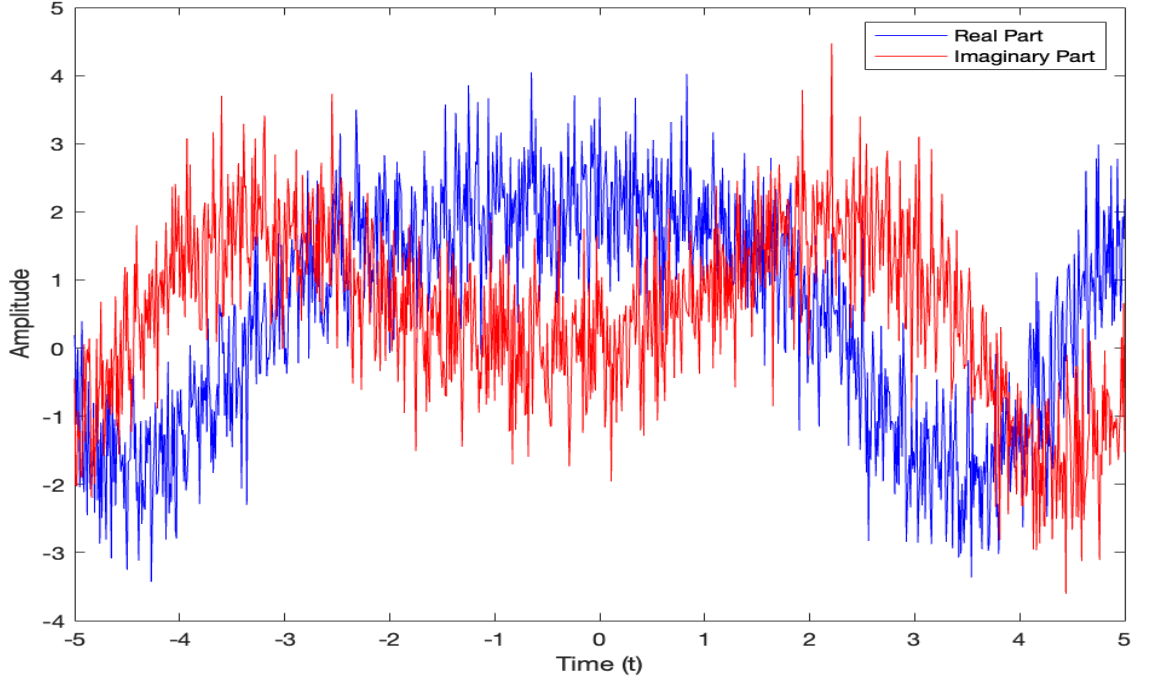

FIGURE 10. Real and Imaginary parts of bi-component signal $u(t) = e^{i(0.1t+0.2t^2)} + e^{i(0.3t+0.2t^2)}$, at 5dB SNR

Now, the AQWD of the auto-terms is given by

$$\mathcal{QW}^{\Lambda}_{f_k}(t,\nu)$$
$$= \begin{cases} |A_k|^2 T|B|e^{i[(D-E)(\nu-t)]} \operatorname{sinc}\left\{\frac{T}{2}\left[(2\xi_k+2A)t+B\nu+\nu_k+\frac{D+E}{2}\right]\right\}, & A=C \\ |A_k|^2|B|e^{i[(A-C)(\nu^2-t^2)+(D-E)(\nu-t)]}\int_{-\frac{T}{2}}^{\frac{T}{2}} e^{i\frac{(C-A)}{4}\Upsilon^2} e^{i\left[(2\xi_k+C+A)t+B\nu+\nu_k+\frac{D+E}{2}\right]\Upsilon} d\Upsilon, & A \neq C. \end{cases}$$

(4.4)

Furthermore

$$f_{k_1}\left(t+\frac{\Upsilon}{2}\right) f_{k_2}^*\left(t-\frac{\Upsilon}{2}\right)$$
$$= A_{k_1} e^{i\left[\nu_{k_1}\left(t+\frac{\Upsilon}{2}\right)+\xi_{k_1}\left(t+\frac{\Upsilon}{2}\right)^2\right]} A_{k_1}^* e^{-i\left[\nu_{k_2}\left(t-\frac{\Upsilon}{2}\right)+\xi_{k_2}\left(t-\frac{\Upsilon}{2}\right)^2\right]}$$
$$= A_{k_1} A_{k_2}^* e^{i\left[\frac{(\xi_{k_2}-\xi_{k_2})}{4}\Upsilon^2+\frac{(\nu_{k_1}+\nu_{k_2})}{2}\Upsilon\right]} e^{i\left[(\xi_{k_2}-\xi_{k_2})t^2+(\nu_{k_1}-\nu_{k_2})t\right]} e^{i(\xi_{k_2}+\xi_{k_2})t\Upsilon}.$$

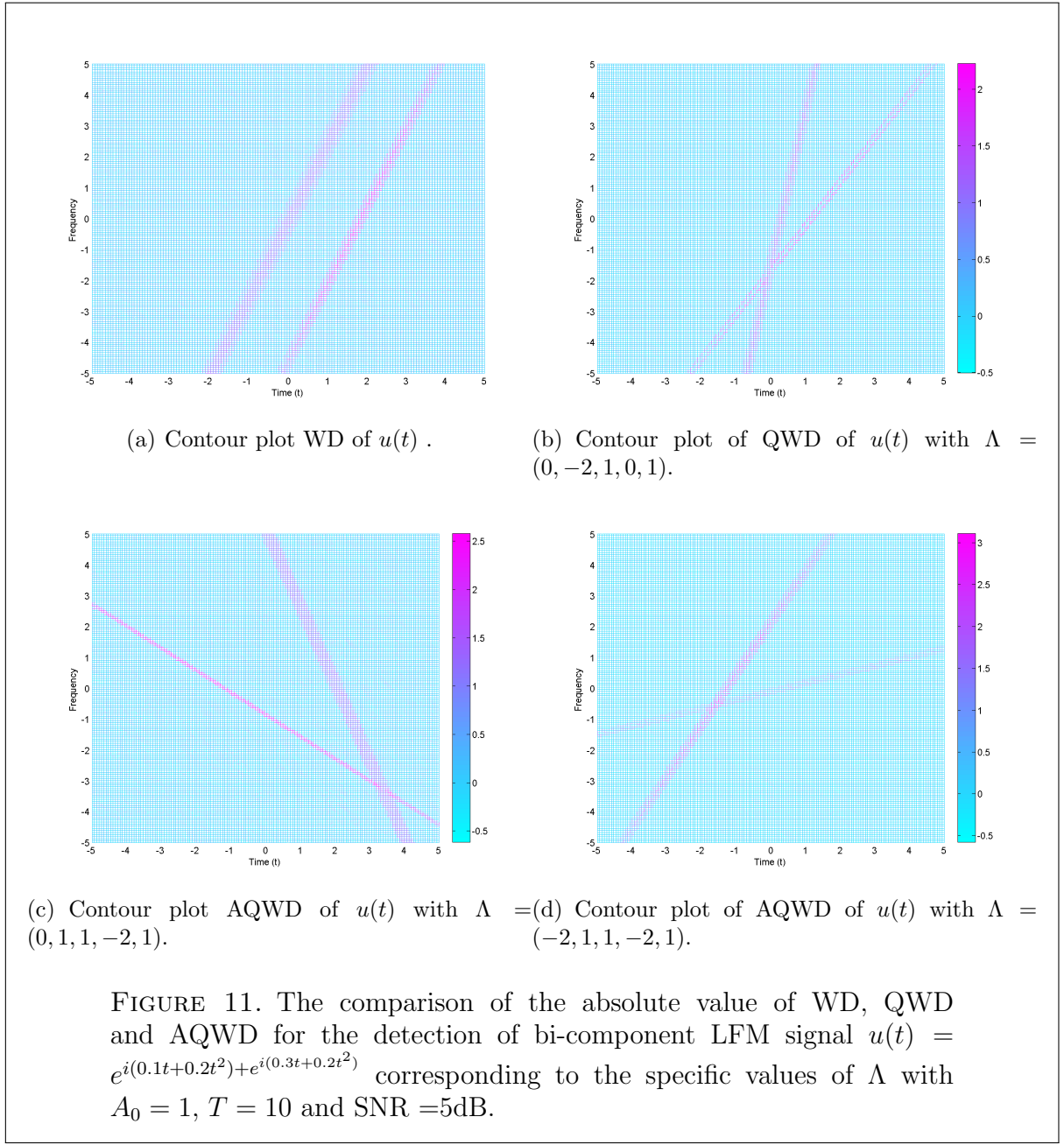

FIGURE 11. The comparison of the absolute value of WD, QWD and AQWD for the detection of bi-component LFM signal $u(t) = e^{i(0.1t+0.2t^2)}+e^{i(0.3t+0.2t^2)}$ corresponding to the specific values of $\Lambda$ with $A_0 = 1$, $T = 10$ and SNR =5dB.

(a) Contour plot WD of $u(t)$.

(b) Contour plot of QWD of $u(t)$ with $\Lambda = (0,-2,1,0,1)$.

(c) Contour plot AQWD of $u(t)$ with $\Lambda = (0,1,1,-2,1)$.

(d) Contour plot of AQWD of $u(t)$ with $\Lambda = (-2,1,1,-2,1)$.

Also, the AQWD of cross-term $\mathcal{QW}^\Lambda_{f_{k_1},f_{k_2}}(t,\nu)$ can be presented as

$$\mathcal{QW}^\Lambda_{f_{k_1},f_{k_2}}(t,\nu) = |B|e^{i\left[(A-C)\nu^2+(D-E)\nu\right]} \int_{-\frac{T}{2}}^{\frac{T}{2}} v_{k_1}\left(t+\frac{\Upsilon}{2}\right) v_{k_2}^*\left(t-\frac{\Upsilon}{2}\right) \times$$
$$e^{i\left[(C-A)\left(t-\frac{\Upsilon}{2}\right)^2+(E-D)\left(t-\frac{\Upsilon}{2}\right)\right]} e^{i[B\nu\Upsilon+2Ct\Upsilon+E\Upsilon]} \mathrm{d}\Upsilon$$
$$= A_{k_1} A_{k_2}^* |B| e^{i\left[(A-C)(\nu^2-t^2)+(D-E)(\nu-t)\right]} e^{i\left[(\xi_{k_2}-\xi_{k_2})t^2+(\nu_{k_1}-\nu_{k_2})t\right]} \int_{-\frac{T}{2}}^{\frac{T}{2}} e^{i\left[\frac{\xi_{k_2}-\xi_{k_2}}{4}+\frac{C-A}{4}\right]\Upsilon^2} \times$$
$$e^{i\left[\left(\xi_{k_2}+\xi_{k_2}+(A+C)\right)t+B\nu+\frac{(D+E)}{2}+\frac{\nu_{k_1}+\nu_{k_2}}{2}\right]\Upsilon} \mathrm{d}\Upsilon$$

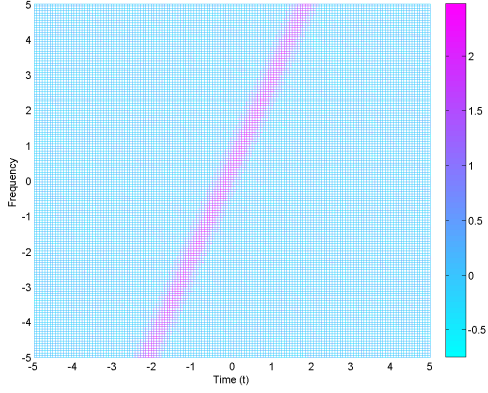
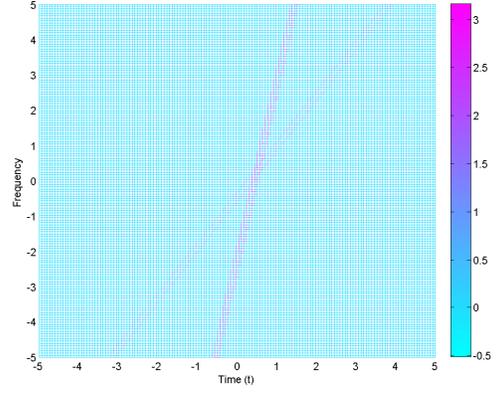

(a) Contour plot AF of $u(t)$.

(b) Contour plot of QAF of $u(t)$ with $\Lambda = (0, -2, 1, 0, 1)$.

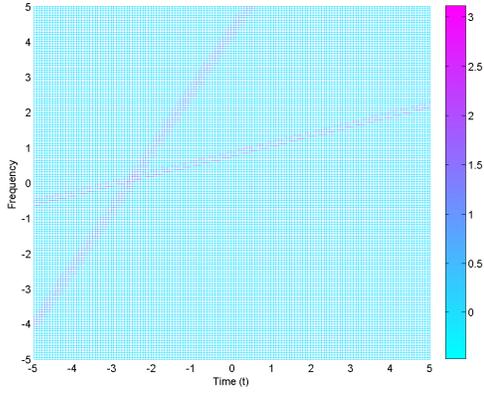
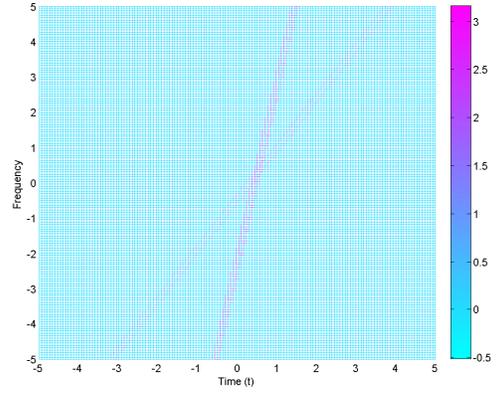

(c) Contour plot of AQAF of $u(t)$ with $\Lambda = (-2, 1, 1, -2, 1)$.

(d) Contour plot of AQAF of $u(t)$ with $\Lambda = (0, -2, 1, 0, 1)$.

FIGURE 12. The comparison of the absolute value of AF, QAF and AQAF for the detection of bi-component LFM signal $u(t) = e^{i(0.1t+0.2t^2)} + e^{i(0.3t+0.2t^2)}$ corresponding to specific values of $\Lambda$ with $A_0 = 1$, $T = 10$ and SNR $=5$dB.

$$= \begin{cases} \mathcal{R}_2(t,\nu) T sinc\left\{\frac{T}{2}\left[(\xi_{k_2}+\xi_{k_2}+A+C)t+B\nu+\frac{\nu_{k_1}+\nu_{k_2}}{2}+\frac{(D+E)}{2}\right]\right\}, & M=0 \\ \mathcal{R}_2(t,\nu) \int_{-\frac{T}{2}}^{\frac{T}{2}} e^{i\left[\frac{\xi_{k_2}-\xi_{k_2}}{4}+\frac{C-A}{4}\right]\Upsilon^2} e^{i\left[(\xi_{k_2}+\xi_{k_2}+(A+C))t+B\nu+\frac{(D+E)}{2}+\frac{\nu_{k_1}+\nu_{k_2}}{2}\right]\Upsilon} d\Upsilon, & M\neq 0, \end{cases}$$
(4.5)

where $M = \frac{\xi_{k_2}-\xi_{k_2}}{4} + \frac{C-A}{4}$ and

$$\mathcal{R}_2(t,\nu) = A_{k_1} A_{k_2}^* |B| e^{i\left[(A-C)(\nu^2-t^2)+(D-E)(\nu-t)\right]} e^{i\left[(\xi_{k_2}-\xi_{k_2})t^2+(\nu_{k_1}-\nu_{k_2})t\right]}.$$

Therefore, it is evident from (4.5) that the AQWD of a cross-term can produce impulses along a straight line $\left[(\xi_{k_2}+\xi_{k_2}+A+C)t+B\nu+\frac{\nu_{k_1}+\nu_{k_2}}{2}+\frac{(D+E)}{2}\right] = 0$ in the $(t,\nu)$-plane when $M = 0$. Then, using (4.4) and (4.5), we conclude that even in the presence of cross-terms, the AQWD still has an influence on the detection performance when $A = C$. The multi-component LFM signal can still be recognized using equation (4.5) even when $A, C$, and $M = 0$. Similarly

$$\mathcal{QA}_f^\Lambda(\Upsilon, \nu) = \sum_{k=1}^n \mathcal{QA}_{f_k}^\Lambda(\Upsilon, \nu) + \sum_{k_1 \neq k_2 = 1}^n \mathcal{QA}_{f_{k_1}, f_{k_2}}^\Lambda(\Upsilon, \nu).$$

The AQAF of auto and cross-terms can be given by

$$\mathcal{QA}_{f_k}^\Lambda(\Upsilon, \nu) = \begin{cases} |A_k|^2 |B| e^{i[(D-E)(\nu - \Upsilon/2)]} e^{i\nu_k \Upsilon} T sinc\left\{\frac{T}{2}\left[(2\xi_k + 2A)\Upsilon + B\nu + E - D\right]\right\}, & A = C \\ \mathcal{R}_3(\Upsilon, \nu) \int_{-\frac{T}{2}}^{\frac{T}{2}} e^{i(C-A)t^2} e^{i[(2\xi_k + C + A)\Upsilon + B\nu + E - D]t} dt, & A \neq C. \end{cases} \quad (4.6)$$

where $\mathcal{R}_3(\Upsilon, \nu) = |A_k|^2 |B| e^{i\left[(A-C)(\nu^2 - \Upsilon^2/4) + (D-E)(\nu - \Upsilon/2)\right]} e^{i\nu_k \Upsilon}$;
and

$$\mathcal{QA}_{f_{k_1}, f_{k_2}}^\Lambda(\Upsilon, \nu)$$
$$= \begin{cases} \mathcal{R}_4(\Upsilon, \nu) T sinc\left\{\frac{T}{2}\left[(\xi_{k_2} + \xi_{k_2} + A + C)\Upsilon + B\nu + \nu_{k_1} - \nu_{k_2} + E - D\right]\right\}, & N = 0 \\ \mathcal{R}_4(\Upsilon, \nu) \int_{-\frac{T}{2}}^{\frac{T}{2}} e^{i[\xi_{k_2} - \xi_{k_2} + C - A]t^2} e^{i[(\xi_{k_2} + \xi_{k_2} + C + A)\Upsilon + B\nu + \nu_{k_1} - \nu_{k_2} + E - D]t} dt, & N \neq 0, \end{cases} \quad (4.7)$$

where $N = \xi_{k_2} - \xi_{k_2} + C - A$ and

$$\mathcal{R}_4(\Upsilon, \nu) = A_{k_1} A_{k_2}^* |B| e^{i\left[(A-C)(\nu^2 - \Upsilon^2/4) + (D-E)(\nu - \Upsilon/2)\right]} e^{i\left[\frac{(\xi_{k_2} - \xi_{k_2})}{4} \Upsilon^2 + \frac{(\nu_{k_1} - \nu_{k_2})}{2} \Upsilon\right]}.$$

According to the relations (4.6) and (4.7), the AQAF is a useful instrument for identifying multi-component LFM signals.

For illustration, the comparison of detecting bi-component LFM signals $u(t) = e^{i(0.1t + 0.2t^2)} + e^{i(0.3t + 0.2t^2)}$ ($|t| \leq 10$) with SNR = 5dB by using Contour plots of WD, QWD and AQWD are displayed in Fig. 11 and Contour plots of AF, QAF and AQAF are displayed in Fig. 12, respectively. As can be seen from Fig. 11 and Fig. 12, the detection performance of AQWD and AQAF for the bi-component LFM signal is significantly more effective.

## 5. Conclusion

This work proposed two advanced time-frequency analysis tools, AQWD and AQAF, as extensions of the classical WD and AF. Key characteristics of these transforms were analyzed, and their effectiveness in identifying single and multi-component LFM signals was demonstrated. Simulation results confirm that AQWD and AQAF offer practical and efficient solutions for LFM signal detection. Moreover, they exhibit greater adaptability than QWD and QAF while delivering superior detection performance compared to WD and AF, making them valuable tools for advanced signal processing applications.

## Declarations

- Ethical Approval: Not Applicable.
- Competing interests: The author has no competing interests.
- Declaration on the Use of AI technologies: The authors affirm that they did not create this article using Artificial Intelligence (AI) technologies.

- Authors' contributions: Not Applicable.
- Funding: No funding was received for this work.